\documentclass[11pt, a4paper, reqno]{amsart}
\usepackage[top=35mm, bottom=35mm, left=30mm, right=30mm]{geometry}
\usepackage{mathptmx}
\usepackage{mathrsfs}
\usepackage{amscd}
\usepackage{verbatim}
\usepackage{url}
\usepackage[all]{xy}
\usepackage{color}
\usepackage[colorlinks=true,citecolor=blue, backref]{hyperref}

\theoremstyle{plain}

\newtheorem{thm}{Theorem}[section]
\newtheorem{lem}[thm]{Lemma}

\newtheorem{cor}[thm]{Corollary}

\theoremstyle{definition}
\newtheorem{defn}[thm]{Definition}

\newtheorem{ques}{Question}[section]
\newtheorem{conj}{Conjecture}

\newtheorem{thma}{Theorem}

\numberwithin{equation}{section}

\newcommand{\Z}{{\mathbb{Z}}}
\newcommand{\N}{\mathbb{N}}
\newcommand{\F}{\mathcal F}

\def \id {{\rm id}}

\begin{document}
\title{Orbits along primes in totally minimal systems}

\author[W.~Huang]{Wen Huang}
\address[W. Huang]{School of Mathematical Sciences, University of Science and Technology of China, Hefei, Anhui 230026, China}
\email{wenh@mail.ustc.edu.cn}

\author[S.~Shao]{Song Shao}
\address[S. Shao]{School of Mathematical Sciences, University of Science and Technology of China, Hefei, Anhui 230026, China}
\email{songshao@ustc.edu.cn}

\author[H.~Xu]{Hui Xu}

\address[H. Xu]{Department of Mathematics, Shanghai Normal University, Shanghai, 200234, China}
\email{huixu@shnu.edu.cn}

\author[X.~Ye]{Xiangdong Ye}
\address[X. Ye]{School of Mathematical Sciences, University of Science and Technology of China, Hefei, Anhui 230026, China}
\email{yexd@ustc.edu.cn}

\thanks{This research is supported by National Key R$\&$D Program of China (No. 2024YFA1013601, 2024YFA1013600), and National Natural Science Foundation of China (12371196, 12426201).}

\maketitle

\begin{abstract}
In this paper, we show that every totally minimal system admits a point whose orbit along the primes is dense. Moreover, assuming the Dickson–Hardy–Littlewood conjecture, we show that in any totally minimal distal system, the orbit along the primes of every point is dense.
\end{abstract}

\section{Introduction}
\subsection{Backgrounds and main results}
By a {\em topological dynamical system}, or simply a {\em system}, we mean a pair $(X,T)$, where $X$ is a compact metric space and $T: X\rightarrow X$ is a homeomorphism. A system $(X,T)$ is transitive if there exists some point $x\in X$ whose orbit $orb(x,T):=\{T^{n}x: n\in \mathbb{Z}\}$ is dense in $X$; it is minimal if the orbit of every point is dense. Furthermore, $(X,T)$ is totally minimal if $(X,T^n)$ is minimal for every positive integer $n$. 

In the study of dynamical systems, we are interested in the question of along which sequences the orbit is equidistributed or dense in the the whole space. 
Let $\mathbb{P}$ denote the set of prime numbers. A system is $\mathbb{P}$-transitive if there is some point whose orbit along the primes is dense. A system is $\mathbb{P}$-minimal if the orbit of every point along the primes is dense. 

It is not surprising that the $\mathbb{P}$-transitivity is deeply related to the distribution of the primes. The $\mathbb{P}$-transitivity of every strongly mixing system is equivalent to the infinitude of the primes, while the distribution of prime orbits in a periodic system is revealed by Dirichlet's Theorem on  the distribution of the primes in arithmetic progressions. For an irrational rotation, Vinogrodov's estimation on exponential sum with prime variables implies that the prime orbits are equidistributed \cite{V37}. Furthermore, Green and Tao shows that this also holds for every minimal connected nilsystem \cite{GT 12}. If a minimal system $(X,T)$ admits a totally ergodic measure, then Bourgain's pointwise convergence theorem \cite{Bour 89} along  the primes shows that $(X,T)$ has dense prime orbits. Unfortunately, not all totally minimal systems admit totally ergodic measures, such as some skew product on tori. Sarnak and Ubis study the prime orbits in horocycle flows  \cite{SU 15}. Kanigowski, Lema\'{n}czky and Radziwill study the prime orbits for analytic skew product on tori \cite{KLR 24}. 


In this paper, we focus on studying the orbits along the primes in a totally minimal system.
In \cite{GHSWY 25}, the authors show that for a totaly minimal system and any integer valued polynomial with degree 2 there is some point whose orbit along the polynomial time being dense. Later this result is proved to be hold for any  integer valued polynomial in \cite{Qiu 23}. In \cite{GHSWY 25}, the authors asked whether every totally minimal system admit a dense prime orbit \cite[Conjecture 2]{GHSWY 25}. 

In this paper, we give an affirmative answer to this question.  

\begin{thma}\label{mainA}
Every totally minimal system is $\mathbb{P}$-transitive.
\end{thma}

For an equicontinuous system, $\mathbb{P}$-transitivity is equivalent to $\mathbb{P}$-minimality. For totally minimal nilsystems, Green and Tao actually show that they are $\mathbb{P}$-minimal \cite{GT 12}. In general, one can not expect that every totally minimal system is $\mathbb{P}$-minimal. Pavlov shows in \cite{Pavlov 08} that for every sequence $A\subset \mathbb{N}$ with zero upper Banach density there is a minimal weakly mixing system $(X,T)$ such that there are uncountably many points $x$'s in $X$ such that $\{T^{n}x: n\in A\}$ is not dense in $X$. Particularly, this holds for the primes. 

Assuming the Dickson-Hardy-Littlewood conjecture holds, we will show that every  totally minimal point-distal system is $\mathbb{P}$-minimal. In particular, this holds for totally minimal distal systems. A finite subset $H\subset \mathbb{Z}$ is {\em admissible} if for every prime $p$, the residual classes of $H$ modulo $p$ does not occupy all residual classes. The {\em Dickson-Hardy-Littlewood conjecture} says that every admissible set can be translated into the primes simultaneously infinitely often. Recall that a system is point-distal if every point is not proximal to other point except itself. 

\begin{thma}\label{mainB}
Assuming the Dickson-Hardy-Littlewood conjecture holds, every totally minimal point-distal system is $\mathbb{P}$-minimal. 
\end{thma}

In fact, under the hypothesis in Theorem \ref{mainB} we will show that for a finite admissible set $A$ and a totally minimal point-distal system $(X,T)$, 
\[ \overline{\{T^{n}x: n\in\mathbb{N},\ A+n \subset \mathbb{P}\}}=X,\ \forall x\in X.\]
In particular, the orbit along the twin primes of every point is dense (see Corollary \ref{cor-3.7}). 


A natural question is to prove Theorem \ref{mainB} without  Dickson-Hardy-Littlewood conjecture.

\begin{conj}
Let $(X,T)$ be a totally minimal point-distal system. Then it is $\mathbb{P}$-minimal, i.e., for each $x\in X$, $\overline{\{T^px:p\in \mathbb{P}\}}=X$.
\end{conj}

Also we have the following question: 

\begin{ques}
Let $(X,T)$ be a totally minimal point-distal system (or distal system) and a non-constant polynomial $p\in \mathbb{Z}[n]$. For each $x\in X$, does $\overline{\{T^{p(n)}x: n\in \Z\}}=X$ hold?
\end{ques}

\medskip

Finally, we remark that the conclusions in this paper also hold if the transformation is a surjective map. In this case, we get an invertible transformation by passing to the natural extension. While both the $\mathbb{P}$-transitivity and $\mathbb{P}$-minimality are preserved under taking factors and inverse limits. 

\subsection{Ideas of the proofs}

Our initial confidence on tackling the question of prime orbits is the observation that Theorem \ref{mainA} is true assuming the Dickson-Hardy-Littlewood conjecture holds. We give this short proof in Section \ref{cond P-tran}. 
 
An important breakthrough on the Dickson-Hardy-Littlewood conjecture starts with Zhang's progress on the twin primes conjecture \cite{Zhang 14} and then Maynard treat arbitrary length of admissible tuples \cite{Maynard 15}. Maynard's result says that prime-producing tuples of any length can be found in any sufficient long admissible tuples. This motivates us to treat systems with enough randomness saying minimal weakly mixing systems (See Section \ref{MT-WM}). 

 
Recently, the third and fourth author introduce the notion of quasi-disjointness and show that every minimal PI system is quasi-disjoint from all minimal systems \cite{XY 27}. This help us pass the $\mathbb{P}$-transitivity of a minimal PI system to the $\mathbb{P}$-transitivity of its maximal equicontinuous factor. While the latter one is a classical known result.

After establishing the $\mathbb{P}$-transitivity of minimal PI systems, it remains to deal with RIC weakly mixing extensions. The method is very similar to the case of minimal weakly mixing systems.  Rather than using Maynard's result directly, we invoke a recent generalization by Bienvenu, Shao and Ter\"{a}v\"{a}inen \cite{BST23}, where they use nilsequences to observe the prime tuples. 
 
 
\section{Preliminaries}
 
In this paper, by a {\em topological dynamical system} or simply a {\em system}, we mean a pair $(X, T)$ with $X$ a compact metric space and $T$ a homeomorphism on $X$. We alway use $\rho$ to denote the metric on the phase space $X$. Sometimes we also consider systems by group actions. Let $X$ be a compact metric space and $G$ be a topological group. Then a continuous $G$-action on $X$ means that there is a continuous map $\phi: G\times X\rightarrow X$ such that $\phi(g, \phi(h,x))=\phi(gh, x)$ and $\phi(e_{G},x)=x$ for any $g,h\in G$ and $x\in X$, where $e_{G}$ is the unit of $G$. For brevity, we write $gx$ for $\phi(g,x)$. 

 We use $\mathbb{Z}$ to denote the set of integers, $\mathbb{N}$ to denote the set of positive integers and $\mathbb{P}$ to denote the set of prime numbers.
 
 \subsection{Basics on topological dynamics}
 
A system $(X,T)$ is {\it transitive} if for any nonempty open subsets $U,V$ of $X$, there is some $n\in\mathbb{Z}$ with $U\cap T^{-n} V\neq\emptyset$. This is equivalent to the existence of transitive point, the orbit of which is dense.  A system $(X,T)$ is {\it minimal} if there is no proper nonempty invariant closed subset, which is equivalent to that every point in $X$ is a transitive point. $(X,T)$ is {\it totally transitive} (resp. {\it totally minimal}) if $(X,T^{n})$ is transitive (resp. minimal) for every $n\in\mathbb{N}$. A system $(X,T)$ is {\it weakly mixing} if the product system $(X\times X, T\times T)$ is transitive. All weakly mixing systems are totally transitive (see \cite[Theorem 1]{CKN 07} for example). 
 
Let $(X,T)$ be a system. Two points $x,y\in X$ are {\it proximal} if $\inf_{n\in\mathbb{Z}}\rho(T^{n}x, T^{n}y)=0$, otherwise they are {\it distal}. $(X,T)$ is distal if  every distinct pair $x\neq y\in X$ are not proximal. $(X,T)$ is {\it equicontinuous} if for every $\epsilon>0$ there is some $\delta>0$ such that $\rho(T^{n}x,T^{n}y)<\epsilon$ for each $n\in\mathbb{Z}$ whenever $\rho(x,y)<\delta$. A point $x\in X$ is a {\it distal point} if $x,y$ are not proximal for every $y\in X$ with $x\neq y$. A system is {\it point-distal} if there is some distal point. Clearly, equicontinuous systems are distal and distal systems are point-distal.

Let $(X,T)$ and $(Y,S)$ be two systems. A {\it homomorphism} or a {\em factor map} between $(X,T)$ and $(Y,S)$ is a continuous surjective map $\pi: X\rightarrow Y$ such that $\pi\circ T=S\circ \pi$. In this case, $(X,T)$ is called an {\it extension} of $(Y,S)$ and $(Y,S)$ is a {\it factor} of $(X,T)$. 
 When $(X,T)$ is minimal, then $\pi$ is semi-open.

 An important fact is that each minimal t.d.s. admits a maximal equicontinuous factor $(X_{eq}, T)$. By Halmos-von Neumann Theorem, a system $(X,T)$ is {minimal and equicontinuous} if and only if $(X,T)$ is topologically conjugate to a {minimal rotation on a compact abelian metric group}. When $(X,T)$ is totally minimal, $X_{eq}$ is also connected.  

  For a subset $A$ of $\mathbb{Z}$ and $x\in X$, we set 
\[ orb_{A}(x):=\{ T^{n}x: n\in A\}.\]

The following result is well known when $A=\mathbb{Z}$. 
\begin{lem}\label{transitive}
Let $(X,T)$ be a system and $A$ be a subset of $\mathbb{Z}$. Then the following statements are equivalent. 
\begin{enumerate}
\item[(1)] There is some point $x\in X$ such that $\overline{orb_{A}(x)}=X$.
\item[(2)] For any nonempty sets $U, V\subset X$, $(\bigcup_{n\in A}T^{-n}U)\cap V\neq\emptyset$.
\item[(3)] There is a dense $G_{\delta}$ subset $\Omega\subset X$ such that for each $x\in \Omega$, $\overline{orb_{A}(x)}=X$.
\end{enumerate}
\end{lem}
\begin{proof}
(1) $\Longrightarrow$ (2). Suppose that $\overline{orb_{A}(x)}=X$. Then for any $k\in \mathbb{Z}$, 
\begin{equation}\label{eq-tran-1}
 \overline{orb_{A}(T^{k}x)}=\overline{T^{k}(orb_{A}(x))}=T^{k}(\overline{orb_{A}(x)})=X.
 \end{equation}
Now let  $U, V\subset X$ be nonempty open sets. Since $\overline{orb_{A}(x)}=X$, there is some $k\in A$ such that $T^{k}x\in V$. By  (\ref{eq-tran-1}), $ \overline{orb_{A}(T^{k}x)}=X$. Then there is some $a\in A$  such that $T^{a}T^{k}x\in U$ and hence $T^{k}x\in (T^{-a}U)\cap V$. Thus we have 
\[ (\bigcup_{n\in A}T^{-n}U)\cap V\neq\emptyset.\]

(2) $\Longrightarrow$ (3). 
Let $\{U_n\}_{n=1}^{\infty}$ be a collection of countable basis for $X$. Set 
\[ \Omega=\bigcap_{n=1}^{\infty} \bigcup_{k\in A} T^{-k}U_{n}.\]
It follows from (2) that $\Omega$ is a dense $G_{\delta}$ set. It is easy to see that for each $x\in\Omega$,  $\overline{orb_{A}(x)}=X$.

Clearly, (3)  implies (1).   
\end{proof}

For $A\subset\mathbb{Z}$, a system is $A$-\emph{transitive} if there is some point whose orbit along $A$ is dense. A system is $A$-\emph{minimal} if there the orbit of every point along $A$ is dense.

\begin{lem}\label{char of A-tran}
Let $(X,T)$ be a system and $A$ be a subset of $\mathbb{Z}$. Then the following two assertions are equivalent.
\begin{enumerate}
\item $(X,T)$ is $A$-transitive.
\item $\overline{\{(x,T^{n}x): x\in X, n\in A\}}=X\times X$.
\end{enumerate}
\end{lem}
\begin{proof}
(1) $\Rightarrow$ (2). Suppose that $(X,T)$ is $A$-transitive. By Lemma \ref{transitive}, there is a dense $G_{\delta}$ subset $\Omega\subset X$ such that $\{T^{n}x: n\in A\}$ is dense in $X$ for each $x\in\Omega$. This implies that $\overline{\{(x,T^{n}x): x\in X, n\in A\}} = X\times X$.

(2) $\Rightarrow$ (1). Suppose that $\overline{\{(x,T^{n}x): x\in X, n\in A\}}=X\times X$. Let $U, V$ be nonempty open subsets of $X$. Then 
\[ (U\times V)\cap \{(x,T^{n}x): x\in X, n\in A\}\neq\emptyset.\]
Thus there is some $x\in X$ and $n\in A$ such that $x\in U$ and $T^{n}x\in V$. Then  $x\in U\cap T^{-n}V$. In particular, $U\cap T^{-n}V\neq\emptyset$. By Lemma \ref{transitive}, $(X,T)$ is $A$-transitive. 
\end{proof}

\subsection{Structure of minimal systems}
Let $\pi: X\rightarrow Y$ be a homomorphism between two systems $(X,T)$ and $(Y,S)$.  For every $k\geq 2$, let 
\[R_{\pi}^{k}:=\{(x_1,\ldots,x_k)\in X^{k}: \pi(x_1)=\cdots=\pi(x_k)\}\]
and $R_{\pi}:=R_{\pi}^{2}$. We call that $\pi$ is 
\begin{enumerate}
\item a {\it proximal extension} if $(x_1,x_2)$ is proximal for every $(x_1,x_2)\in R_{\pi}$;
\item an {\it equicontinuous extension} if for every $\epsilon>0$ there is some $\delta>0$ such that $\rho(T^{n}x_1,T^{n}x_2)<\epsilon$ for each $n\in\mathbb{Z}$ and every  $(x_1,x_2)\in  R_{\pi}$ with $\rho(x,y)<\delta$;
\item an {\it almost one to one extension} if the set $\{x\in X: |\pi^{-1}\pi(x)|=1\}$ is a dense $G_{\delta}$ subset of $X$;
\item a {\it weakly mixing extension} if the subsystem $(R_{\pi}, T\times T)$ of $(X\times X, T\times T)$ is transitive;
\item a {\it RIC extension} if $\pi$ is open and $R^{k}_{\pi}$ has dense minimal points for every $k\geq 2$.  
\end{enumerate}
The definition of RIC extension above is an equivalent form (\cite[Theorem 1.9]{AMWW 84}) rather that the original form that defined by Ellis semigroup (\cite{EGS 75}). The following result is well known to experts. As we can not find a reference, we prove it for completeness.

\begin{lem}\label{total min for prox ext}
Let $\pi: (X,T)\rightarrow (Y,S)$ be a proximal extension between minimal systems. If $(Y,S)$ is totally minimal then $(X,T)$ is also totally minimal.
\end{lem} 
 \begin{proof}
 We first show that $\pi: (X,T^{k})\rightarrow (Y,S^{k})$ is also a proximal extension for each $k\in\mathbb{N}$. To this end, we fix $(x_1,x_2)\in R_{\pi}$ and $k\in \mathbb{N}$. We need to show that $x_1$ and $x_2$ is a proximal under $T^{k}$. Further, we need to show that for every $\epsilon>0$, there is some $m\in\mathbb{Z}$ such that $\rho(T^{kn}x_1, T^{kn}x_2)<\epsilon$. Since $x_1$ and $x_2$ is a proximal under $T$, the set $\{n\in\mathbb{Z}: \rho(T^{n}x_1, T^{n}x_2)<\epsilon\}$ is a thick subset of $\mathbb{Z}$ (see \cite[Lemma 8.1]{Fur81}), which means it contains intervals of arbitrary length. Thus there is some $n \in\mathbb{Z}$ such that $\rho(T^{kn}x_1, T^{kn}x_2)<\epsilon$. 
 
 On the one hand, the minimality of $(X,T)$ implies that $X$ is the union of finitely many minimal sets of $T^{k}$ for every $k$.  On the other hand, $(X, T^{k})$ has a unique minimal set since $\pi: (X,T^{k})\rightarrow (Y,S^{k})$ is also a proximal extension (see \cite[Lemma 2.1 ]{G75}). These two aspects imply that  $(X, T^{k})$ is minimal for each $k\in\mathbb{N}$. Thus $(X,T)$ is  totally minimal.
 \end{proof}
 
The following result was first gotten by McMahon \cite[Theorem 1.3]{Mc76}, and see also \cite[Theorem 2.12]{AMWW 84}, \cite[Theorem 2.7]{G05} or \cite[Theorem 2.6.4]{V77}. 

\begin{thm}
If $\pi: (X,T)\rightarrow (Y,S)$ is a RIC weakly mixing extension between minimal systems, then  $(R_{\pi}^{k}, T\times\cdots\times T)$ is transitive for each $k\geq 2$.
\end{thm}

Let $\pi: (X,T)\rightarrow (Y,S)$ be an extension between minimal systems.  Suppose that there is a countable ordinal $\eta$ and a family of systems $(X_{\alpha}, T_{\alpha})_{\alpha\leq \eta}$ and homomorphisms $\pi_{\alpha,\beta}: X_{\alpha}\rightarrow X_{\beta}, \beta<\alpha\leq \eta$ such that 
\begin{enumerate}
\item[(a)] $(X,T)=(X_{\eta}, T_{\eta}), \pi=\pi_{\eta,0}$;
\item[(b)] $\pi_{\alpha,\gamma}=\pi_{\beta,\gamma}\pi_{\alpha,\beta}$ for any $\gamma<\beta<\alpha\leq\eta$;
\item[(c)] $(X_{\alpha}, T_{\alpha})=\underset{\longleftarrow}{\lim}_{\beta<\alpha}(X_{\beta}, T_{\beta})$ when $\alpha$ is a limit ordinal. 
\end{enumerate}
 We say that $\pi$ is 
 \begin{enumerate}
 \item[(i)] an {\bf I}-{\it extension} if $\pi_{\alpha+1,\alpha}$ is an equicontinuous extension for each $\alpha<\eta$;
 \item[(ii)] an strict {\bf PI}-{\it extension}  if $\pi_{\alpha+1,\alpha}$ is an equicontinuous  or proximal  extension for each $\alpha<\eta$;
  \item[(iii)] a strict {\bf AI}-{\it extension}  if $\pi_{\alpha+1,\alpha}$ is an equicontinuous  or almost one to one extension for each $\alpha<\eta$.
 \end{enumerate}
 
 A minimal system $(X,T)$ is a strict {\bf PI} system if it a strict {\bf PI}-extension of the singleton system; it is a {\bf PI} system if there is a minimal proximal extension $\theta: (\tilde{X}, \tilde{T})\rightarrow (X,T)$ such that $(\tilde{X}, \tilde{T})$ is a strictly {\bf PI} system.  A minimal system $(X,T)$ is a strict {\bf AI} system if it a strictly {\bf AI}-extension of the singleton system; it is an {\bf AI} system if there is an almost one to one extension $\theta: (\tilde{X}, \tilde{T})\rightarrow (X,T)$ such that $(\tilde{X}, \tilde{T})$ is a strict {\bf AI} system.

Furstenberg showed that a minimal system is distal if and only if it is an {\bf I}-extension of a trivial system \cite{Fur 63}. Veech showed that a minimal system is point-distal   if and only if it is an {\bf AI}  system \cite{Veech 70}. 


Finally we have the structure theorem for minimal systems
(see Ellis-Glasner-Shapiro \cite{EGS 75},
McMahon \cite{Mc76}, Veech \cite{V77}).

\begin{thm}[Structure theorem for minimal systems]\label{structure}
Let $(X,T)$ be a minimal system. Then we have the following diagram:
\[
\xymatrix{
  X          & X_\infty \ar[d]^{\pi_\infty} \ar[l]^{\theta} \\
                & Y_\infty             }
\]

\noindent where $X_\infty$ is a proximal extension of $X$ and a RIC
weakly mixing extension of the strictly PI-system $Y_\infty$.

The extension $\pi_\infty$ is an isomorphism (so that
$X_\infty=Y_\infty$) if and only if  $X$ is a PI-system.
\end{thm}

\subsection{Regionally proximal relations}
  Let $(X,T)$ be a  system and $d$ be a positive integer. The points $x,y\in X$ are said to be \emph{regionally proximal of order} $d$ if for any $\delta>0$, there are $x',y'\in X$ and a vector ${\bf n}=(n_1,\ldots,n_d)\in\mathbb{Z}^{d}$ such that $d(x,x')<\delta, d(y,y')<\delta$, and 
  \[d(T^{n_1\epsilon_1+\cdots +n_d\epsilon_d}x',T^{n_1\epsilon_1+\cdots +n_d\epsilon_d}y')<\delta, \ \forall (\epsilon_1,\ldots,\epsilon_d)\in\{0,1\}^{d}\setminus \{(0,0,\ldots,0)\}.\]
 The set  of regionally proximal pairs of order $d$ is denoted by ${\bf RP}^{[d]}(X)$. 
  
  The notion of regionally proximal of order of $d$ is introduced by Host, Kra and Maass \cite{HKM}. They also show that ${\bf RP}^{[d]}(X)$ is an invariant equivalence closed relation when $(X,T)$ is a minimal distal system. Shao and Ye show that ${\bf RP}^{[d]}(X)$ is an invariant equivalence closed relation for any minimal system $(X,T)$ \cite{SY 12}.  By observing that ${\bf RP}^{[1]}(X)\supset {\bf RP}^{[2]}(X)\supset\cdots$, we define 
  \[ {\bf RP}^{[\infty]}(X)=\bigcap_{d=1}^{\infty} {\bf RP}^{[d]}(X).\]
  
  Let $\phi: Y_1\rightarrow Y_2$ be a map between two sets $Y_1$ and $Y_2$. A subset $L\subset Y_1$ is called $\phi$-{\it saturated} if $\phi^{-1}\phi(L)=L$.  We need the following saturation theorem of Glasner, Huang, Shao, Weiss and Ye. 
  \begin{thm}\cite[Theorem 4.2]{GHSWY 25}\label{saturated thm}
  Let $\pi: (X,T)\rightarrow (Y,S)$ be an  extension between minimal systems. If $\pi$ is open and ${\bf RP}^{[\infty]}(X)\supset R_{\pi}=\{(x,x')\subset X\times X: \pi(x)=\pi(x')\}$, then for any $d\in\mathbb{N}$, there is a dense $G_{\delta}$ subset $\Omega$ of $X$ such that the set 
  \[ \overline{\{(T\times T^2\times\cdots\times T^{d})^{n}(x,x,\ldots,x):  n\in\mathbb{N}\}}\]
  is $(\pi\times\pi\times\cdots\times\pi)$-saturated, for any $x\in\Omega$.
  \end{thm}

  \subsection{Dickson-Hardy-Littlewood conjecture}
 
\begin{defn}
A tuple $(h_1,h_2,\ldots,h_k)$ or a  subset $\{h_1,h_2,\ldots,h_k\}$ of integers is {\it admissible} if for any prime $p$, there is some $a$ such that $a\not\equiv h_i\mod p$ for each $i\in\{1,\ldots,k\}$. It is {\it prime-producing} if there are infinitely many $n$ such that $n+h_1,\ldots,n+h_k$ are simultaneously prime. 
\end{defn}

For example,  for every odd prime $p$, there is a non quadratic residue. Thus for any $n$, the tuple $\left(1^{2}, 3^2,\ldots, (2n+1)^{2}\right )$ is admissible. The following fact will be used frequently. 

\begin{lem}\label{W-trick}
Let $k\in \mathbb{N}$ and $\displaystyle w=\prod_{p\in\mathbb{P}, p\leq k}p$. Then $\{wn_1, 2wn_2, \ldots, kwn_k\}$ is admissible for any integers $n_1, n_2,\ldots,n_k$.
\end{lem}

\begin{proof}
Fix a prime $p$ and integers $n_1, n_2,\ldots,n_k$. If $p>k$ then it clear that $\{wn_1, 2wn_2, \ldots, kwn_k\}$ is admissible, since there are $p$ residual classes modulo $p$. If $p\leq k$ then $wn_1, 2wn_2, \ldots, kwn_k$ are multiples of $p$ and hence the residual classes of $\{wn_1, 2wn_2, \ldots, kwn_k\}$ cannot occupy all residual classes modulo $p$. Therefore, $\{wn_1, 2wn_2, \ldots, kwn_k\}$ is admissible.
\end{proof}

\begin{conj}[Dickson-Hardy-Littlewood conjecture]
Every admissible tuple $(h_1,\ldots,h_k)$ is prime-producing.
\end{conj}
This is a special case of Dickson's conjecture \cite{Dick 04} and Hardy and Littlewood conjectured a more quantitative estimate \cite{HL 23}. 

The following result is obtained by Maynard and independently by Tao. 
\begin{thm}[\cite{Maynard 15}]\label{MT}
For any $k\in\N$, there is some $m=m_{k}$ such that for any admissible sequence $(h_1,h_2,\ldots,h_{m})$, there is some $n\in\mathbb{N}$ and $1\leq i_1<\ldots<i_k
\leq m$ such that $\{h_{i_1},\ldots,h_{i_k}\}+n\subset\mathbb{P}$.
\end{thm}
Theorem \ref{MT} helps us establish the $\mathbb{P}$-transitivity of every minimal weakly mixing system. To deal with RIC weakly mixing extensions, we need a generalization of Theorem \ref{MT} obtained by Bienvenu, Shao and Ter\"{a}v\"{a}inen \cite[Proposition 5.1]{BST23}.  Avoiding introducing extra notations, we do not restate their result here but to give an immediate consequence of \cite[Proposition 5.1]{BST23} without proof. 

\begin{lem}\label{torus}
For any $k\in \mathbb{N}$, there is some $m=m_k\in\mathbb{N}$ such that for any $d\in\mathbb{N}$, any minimal rotation $R$ on the torus $\mathbb{T}^{d}$ and any admissible set $H$ with $|H|=m$, one has
\[\overline{\{R^{n}x: n\in\mathbb{N}, |(H+n)\cap\mathbb{P}|\geq k\}}=\mathbb{T}^{d}, \ \forall x\in\mathbb{T}^{d}.\]
\end{lem}

\subsection{Family of subsets of integers}
In this subsection, we introduce the notion of families of subsets of $\N$. All concepts can be defined for $\Z$ in the same way.

A {\it family}  $\mathcal{F}$ of $\mathbb{N}$ is a collection of nonempty subsets of $\mathbb{N}$ that is upwards hereditarily, that is for any $E\subset F\subset \mathbb{N}$, if $E\in \mathcal{F}$ then $F\subset \mathcal{F}$. The dual family of $\mathcal{F}$  is define to be 
\[\mathcal{F}^{*}:=\{ A\subset \mathbb{N}: A\cap F\neq\emptyset, \forall F\in \mathcal{F}\}.\]

A subset $A\subset \mathbb{N}$ is {\it syndetic} if it has a bounded gap, that is there is some $L\in \mathbb{N}$ such that $a_{i+1}-a_i\leq L$ for each $i\in\mathbb{N}$, where $a_1<a_2<a_3<\cdots$ is an ordering of $A$. A subset $A\subset \mathbb{N}$ is {\it thick} if it contains intervals of arbitrary length, that is for any $\ell\in\mathbb{N}$ there is some $n\in\mathbb{N}$ such that $ \{n,n+1,\ldots, n+\ell\}\subset A$. A subset $A\subset \mathbb{N}$ is {\it piecewise syndetic} if there is a syndetic subset $B\subset \mathbb{N}$ and thick subset $C\subset \mathbb{N}$ such that $A=B\cap C$. A subset $A\subset \mathbb{N}$ is {\it thickly syndetic} if  intervals of arbitrary length can occur in $A$ syndetic often, that is for any $\ell\in\mathbb{N}$, there is a syndetic subset $B\subset \mathbb{N}$ such that 
\[\{b,b+1,\ldots, b+\ell\}\subset A, \ \forall b\in B.\] 

 Let $\mathcal{F}_{s}, \mathcal{F}_{t},\mathcal{F}_{ps}, \mathcal{F}_{ts}$ denote the collection of syndetic subsets, thick subsets, piecewise syndetic subsets and  thickly syndetic subsets, respectively. These are all families of $\mathbb{N}$. It is well-known that $\mathcal{F}_{s}$ and  $\mathcal{F}_{t}$ are dual to each other, and $\mathcal{F}_{ps}$ and  $\mathcal{F}_{ts}$ are dual to each other.

A family $\mathcal{F}$ of $\mathbb{N}$ is a {\it filter} if it has {\it finite intersection property}, that is  for any $F_1,F_2\in\mathcal{F}$ there is some $F_3\in\mathcal{F}$ such that $F_3\subset F_1\cap F_2$.
 A family $\mathcal{F}$ of $\mathbb{N}$ has {\it Ramsey property} if $A\cup B\in \mathcal{F}$ then either $A\in\mathcal{F}$ or $B\in\mathcal{F}$. 

\begin{lem}[{\cite[Lemma 9.5]{Fur81}}]\label{filter vs Ramsey}
A family $\mathcal{F}$ of $\mathbb{N}$ is a filter if and only if its dual family has Ramsey property.
\end{lem}

A collection $\mathcal{C}$ of nonempty subsets of $\mathbb{N}$ will generate a family 
\[\mathcal{F}_{\mathcal{C}}=\{E\subset \mathbb{N}: \exists F\in \mathcal{C}, F\subset E\}.\]
Thus if $\mathcal{C}$ has finite intersection property then the family $\mathcal{F}_{\mathcal{C}}$ is a filter. 	In this case, we say that $\mathcal{C}$ generates a filter. 

The following lemma is a consequence of Lemma \ref{filter vs Ramsey}. 
\begin{lem}\label{filter}
Let $\mathcal{C}$ be  a collection of nonempty subsets of  $\mathbb{N}$. If $\mathcal{C}$ generates a filter, then for any finite partition $\mathbb{N}=A_1\cup A_2\cup\cdots\cup A_r$, there is some $i\in\{1,2,\ldots,r\}$ such that 
\[ A_i\cap C\neq\emptyset, \ \forall C\in \mathcal{C}.\]
\end{lem}

\section{Conditional results}

In this section, we establish some conditional results on $\mathbb{P}$-transitivity and $\mathbb{P}$-minimality. Assuming the Dickson-Hardy-Littlewood conjecture holds, we show that every totally minimal system is $\mathbb{P}$-transitive and every totally minimal point-distal system is $\mathbb{P}$-minimal. 

\subsection{$\mathbb{P}$-transitivity of totally minimal systems}\label{cond P-tran}

Let $(X,T)$ be a system and $d\in\mathbb{N}$. Define
\[\sigma_{d}(T)=T\times T\times \cdots\times T\ \ \text{ and }\ \ \tau_{d}(T)=T\times T^{2}\times\cdots\times T^{d}.\]
Further, let 
\[ N_{d}(X, T)=\overline{\{ \tau_{d}(T)^{n}(x,\ldots,x): n\in\mathbb{Z}, x\in X\}}.\]

The following result is due to Glasner \cite{Glasner 94}, and see \cite[Proposition 1.55]{G03} for another proof.

\begin{lem}[{\cite[Theorem 5.1]{Glasner 94}}]\label{Glasner}
If $(X,T)$ is minimal then $(N_{d}(X,T), \langle \sigma_d, \tau_{d}\rangle)$ is also minimal for any $d\in\mathbb{N}$.
\end{lem}

\begin{lem}[{\cite[Theorem C]{GHSWY 25}}]\label{N_d for totally minimal}
Let $(X,T)$ be a minimal system and $k\in\mathbb{N}$. Then $(X,T^{k})$ is minimal if and only if $N_{d}(X,T)=N_{d}(X,T^{k})$ for each $d\in\N$.
\end{lem}

The following observation is followed from the definition.
\begin{lem}\label{observation}
Let $(X,T)$ be a minimal system and $d\in\mathbb{N}$. For any nonempty open set $U\subset X$, one has
\[ N_{d+1}(X,T)\cap (U\times T^{-1}U\times\cdots\times T^{-d}U)\neq\emptyset.\]
\end{lem}

\begin{thm}
Let $(X,T)$ be a totally minimal system. Assuming the Dickson-Hardy-Littlewood conjecture, there is some point $x\in X$ such that $\overline{orb_{\mathbb{P}}(x)}=X$.
\end{thm}
\begin{proof}
By Lemma \ref{transitive}, it suffices to show that for any nonempty open sets $U, V$ of $X$, $(\bigcup_{p\in\mathbb{P}}T^{-p}U)\cap V\neq\emptyset$. 

Since $(X,T)$ is minimal, there is some $d\in \N$ such that 
\begin{equation}\label{eq2}
X=V\cup T^{-1}V\cup \cdots\cup T^{-d}V.
\end{equation}
Let
\[ w= \prod_{p\in\mathbb{P}, p\leq d+1}p.\]
By Lemma \ref{W-trick},   $(0,wm, 2wm, \ldots, dwm)$ is admissible for any $m\in\N$.
By Lemma \ref{observation} and Lemma \ref{N_d for totally minimal}, one has 
\[ N_{d+1}(X, T^{w})\cap (U\times T^{-1} U\times \cdots\times T^{-d}U)\neq\emptyset.\]
Recall that $(N_{d+1}(X,T^{w}), \langle \sigma_{d+1}(T^{w}), \tau_{d+1}(T^{w})\rangle)$ is  minimal. There is some $x\in X$ and $n,m\in\mathbb{N}$ such that 
\begin{equation}\label{eq3}
 T^{wn} x\in  U, T^{wm+wn}x\in T^{-1}U, T^{2wm+wn}x\in T^{-2}U, \ldots, T^{dwm+wn}x\in T^{-d}U.  
 \end{equation}

\noindent{\bf Claim}. For any $k\in\mathbb{Z}$, 
\[(U\cup T^{-wm}U\cup T^{-2wm}U\cup\cdots\cup T^{-dwm}U)\cap T^{-k}V\neq\emptyset.\]

\begin{proof}[Proof of Claim] 
 By (\ref{eq3}), one has 
 \[T^{wn}x\in U, T^{wn+1}x\in T^{-wm}U, T^{wn+2}x\in T^{-2wm}U,\ldots, T^{wn+d}x\in T^{-dwm}U.\]
By (\ref{eq2}), there is some $i\in\{0,1,\ldots,d\}$ such that $T^{wn+k}x\in T^{-i}V$. Thus $T^{wn+i}x\in T^{-k}V$. Therefore, 
\[(U\cup T^{-wm}U\cup T^{-2wm}U\cup\cdots\cup T^{-dwm}U)\cap T^{-k}V\neq\emptyset.\]
Thus we complete the proof of Claim.
\end{proof}

Recall that $(0,wm, 2wm, \ldots, dwm)$ is admissible. By the Dickson-Hardy-Littlewood conjecture, there some some $k\in \mathbb{Z}$ such that $\{k,wm+k,2wm+k,\ldots, dwm+k\}\subset \mathbb{P}$. Then it follows from the Claim that $(\bigcup_{p\in\mathbb{P}}T^{-p}U)\cap V\neq\emptyset$.  This completes the proof.
\end{proof}

\subsection{$\mathbb{P}$-minimality of totally minimal point-distal systems }
In this subsection we show the $\mathbb{P}$-minimality of totally minimal point-distal systems assuming the Dickson-Hardy-Littlewood conjecture.

\begin{defn}
A system $(X,T)$ has {\it admissible finite cover property} (AFC for short) if for every nonempty open set $U\subset X$ and every nonempty finite admissible set $A\subset \mathbb{Z}$, there is a finite set $F\subset \mathbb{Z}_{+}$ such that $A+F$ is admissible and 
\[ X=\bigcup_{n\in F} T^{-n}U.\]
\end{defn}

\begin{lem}\label{AFC}
Every totally minimal point-distal system has AFC. 
\end{lem}

At present, we show Theorem \ref{mainB} assuming Lemma \ref{AFC}, whose proof is given in the end of this subsection.

\begin{proof}[Proof of Theorem \ref{mainB}]
Let $(X,T)$ be a totally minimal point-distal system. To show it is $\mathbb{P}$-minimal, it suffices to prove that for any nonempty open set $U\subset X$ we have $X=\bigcup_{p\in\mathbb{P}}T^{-p}U$.

By Lemma \ref{AFC}, $(X,T)$ has AFC. Since $A=\{0\}$ is admissible, there is a finite subset $F\subset \mathbb{Z}_{+}$ such that $X=\bigcup_{n\in F}T^{-n}U$ and $F=A+F$ is admissible. By Dickson-Hardy-Littlewood conjecture, there is some $k\in\mathbb{Z}_{+}$ such that $F+k$ is contained in the primes. Thus by $X=\bigcup_{n\in F}T^{-n}U$ we have
\[ X=\bigcup_{n\in F+k}T^{-n}U\subset \bigcup_{p\in\mathbb{P}}T^{-p}U.\]
The proof of Theorem \ref{mainB} is complete.
\end{proof}

\begin{cor}\label{cor-3.7}
Let $(X,T)$ be  a totally minimal point-distal system. Assuming the Dickson-Hardy-Littlewood conjecture, for every finite admissible set $A$ and $x\in X$, the set 
\[\{T^{n}x: n+A\subset \mathbb{P}\}\]
is dense in $X$. In particular, the orbit of every point along the twin primes is dense. 
\end{cor}

\begin{proof}
By Lemma \ref{AFC}, $(X,T)$ has AFC. For any nonempty open set $U\subset X$, there is a finite set $F$ such that $A+F$ is adimissible and
\[X=\bigcup_{n\in F} T^{-n}U.\]
Thus there is some $k\in\mathbb{Z}$ such that $A+F+k\subset \mathbb{P}$. Clearly, $F+k\subset \{n\in \mathbb{Z}: n+A\subset \mathbb{P}\}$ and 
\[X=\bigcup_{n\in F} T^{-(n+k)}U\subset \bigcup_{n, n+A\subset\mathbb{P}}T^{-n}U.\]
The arbitrary choice of $U$ shows that $\{T^{n}x: n+A\subset \mathbb{P}\}$ is dense in $X$ for every $x\in X$.
\end{proof}

 In the rest part of this section, we show Lemma \ref{AFC}. It suffices to show that APC is preserved by isometric extensions, almost one to one extensions and inverse limits for totally minimal systems.
 
 \begin{defn}
 An extension $\pi: (X,T)\rightarrow (Y,S)$ is \emph{admissible} if for any nonempty open set $U$ of $X$ and a finite admissible set $A\subset \mathbb{Z}$, there some finite subset $F\subset \mathbb{Z}_{+}$ and a point $x_0\in X$ such that $A+F$ is admissible and 
 \[ \pi^{-1}\pi(x_0)\subset \bigcup_{n\in F} T^{-n} U.\]
 \end{defn}
 \begin{lem}[Admissible extensions]\label{admissible ext}
 Let $\pi: (X,T)\rightarrow (Y,S)$ be an admissible extension between minimal systems. If $(Y,S)$ has AFC then so does $(X,T)$.
 \end{lem}
 \begin{proof}
 Let $U$ be a nonempty open set of $X$ and $A\subset \mathbb{Z}$ be a finite admissible set. Since $\pi$ is admissible, there is a finite subset $F_1\subset \mathbb{Z}_{+}$ and $x_0\in X$ such  that $A+F_1$ is admissible and 
\[\pi^{-1}\pi(x_0)\subset \bigcup_{n\in F_1}T^{-n}U.\] 

Now take a nonempty open  subset  $O$ of $Y$  such that 
\[ \pi^{-1}(O)\subset \bigcup_{n\in F_1}T^{-n}U.\]
Applying the AFC property of $(Y,S)$ with the admissible set $A+F_1$ and open set $O$, there is some finite subset $F_2\subset \mathbb{Z}_{+}$ such that $A+F_1+F_2$ is admissible and 
\[ Y=\bigcup_{n\in F_2}S^{-n}O.\]
Now let $F=F_1+F_2$. It remains to show that 
\[ X=\bigcup_{n\in F}T^{-n}U.\]
For any $x\in X$, there is some $k\in F_2$ such that $\pi(x)\in S^{-k}O$ and hence 
\[ x\in \pi^{-1}(S^{-k}O)=T^{-k}\pi^{-1}(O)\subset T^{-k}\bigcup_{n\in F_1}T^{-n}U\subset \bigcup_{n\in F}T^{-n}U.\]
This shows the AFC property of $(X,T)$.
 \end{proof}

 \begin{lem}[Equicontinuous extensions]\label{EE}
 Let $\pi: (X,T)\rightarrow (Y,S)$ be an isometric extension between totally minimal systems. Then $\pi$ is an admissible extension. 
 \end{lem}
\begin{proof}
Let $U$ be a nonempty open set of $X$ and $A\subset \mathbb{Z}$ be a finite admissible set. Fix $x_0\in X$.
We show that there is a finite subset $F\subset \mathbb{Z}_{+}$ such that $A+F$ is admissible and 
\[\pi^{-1}\pi(x_0)\subset \bigcup_{n\in F}T^{-n}U.\]

Take a nonempty open set $V$ of $X$ and $\epsilon>0$ such that $B_{\epsilon}(V)=\{x\in X: \rho(x,V)<\epsilon\}\subset U$. Since $\pi$ is an isometric extension, there is some $\delta>0$ such that for any $x,x'\in X$ with $\pi(x)=\pi(x')$, one has 
\[d(T^nx, T^nx')<\epsilon \ \text{ for any $n\in\mathbb{Z}$ whenever }\ d(x,x')<\delta.\]
Since $\pi^{-1}\pi(x_0)$ is compact, there are $x_1,x_2,\ldots,x_r\in \pi^{-1}\pi(x_0)$ such that 
\[ \pi^{-1}\pi(x_0)\subset \bigcup_{i=1}^{r}B_{\delta}(x_i).\]
Set 
\[w:=\prod_{p\leq r |A|, p\in\mathbb{P}} p. \]
Since $(X,T)$ is totally minimal, there are $n_1,\ldots, n_r\in \mathbb{N}$ such that $T^{wn_1}x_1,\ldots, T^{wn_r}x_r\in V$. Then we have 
\[\pi^{-1}\pi(x_0)\subset\bigcup_{i=1}^{r}T^{-wn_i}U.\]
 Indeed,  any $x\in \pi^{-1}\pi(x_0)$ lies in $B_{\delta}(x_i)$ for some $i=1,\ldots,r$ and hence  
 \[T^{wn_i}x\in B_{\epsilon}(T^{wn_i}x_i)\subset B_{\epsilon}(V)\subset U.\]
 
  Let $F_1=\{wn_1,\ldots wn_2\}$. It remains to show that $A+F_1$ is admissible. But this is clear, since $|A+F_1|\leq r|A|$ and $A$ is admissible. In fact, for a prime $q$, if $q> r|A|$ then it is clear that $A+F_1$ misses some residual class modulo $q$; if $q\leq r|A|$ then the residual classes of $A+F_1$ modulo $q$ coincide with $A$'s and hence they cannot occupy all residual classes.
\end{proof}

 \begin{lem}[Almost one to one extensions]\label{AE}
 Let $\pi: (X,T)\rightarrow (Y,S)$ be an almost one to one extension between minimal systems. Then $\pi$ is an admissible extension. 
 \end{lem}
 \begin{proof}
Let $U$ be a nonempty open set of $X$ and $A\subset \mathbb{Z}$ be a finite admissible set. Since $\pi$ is almost one to one, there is some $x_0$ with $\pi^{-1}\pi(x_0)=\{x_0\}$. Let $F=\{0\}$. Then $A+F=A$ is admissible and 
\[\pi^{-1}\pi(x_0)=\{x_0\} \subset U=\bigcup_{n\in F}T^{-n}U.\]
 Thus $\pi$ is admissible.
 \end{proof}
 
\begin{lem}[Inverse limits]\label{inv lim}
Let $(X,T)=\varprojlim_{\alpha\in \Lambda}(X_{\alpha}, T_{\alpha})$ be a totally minimal inverse limit system. If each $(X_{\alpha}, T_{\alpha})$ has AFC then so does $(X,T)$. 
\end{lem} 
\begin{proof}
Let $U$ be a nonempty open set of $X$ and $A\subset \mathbb{Z}$ be a finite admissible set. Since the topology on $X$ is determined by $\pi_{\alpha}, \alpha\in \Lambda$, there is some $\alpha\in \Lambda$ and a nonempty open set $U_{\alpha}$ of $X_\alpha$ such that $\pi^{-1}_{\alpha}(U_{\alpha})\subset U$. It follows from the AFC property of $(X_{\alpha},T_{\alpha})$ that there is some finite subset $F\subset \mathbb{Z}_{+}$ such that $A+F$ is admissible and
\[X_{\alpha}=\bigcup_{n\in F}T_{\alpha}^{-n}U_{\alpha}.\]
Thus 
\[ X=\pi_{\alpha}^{-1}(X_\alpha)=\bigcup_{n\in F}\pi_{\alpha}^{-1}(T_{\alpha}^{-n}U_{\alpha})=\bigcup_{n\in F}T^{-n}\pi_{\alpha}^{-1}(U_{\alpha})\subset \bigcup_{n\in F}T^{-n}U\subset X.\]
This shows that $(X,T)$ has AFC.
\end{proof}

Now we can finish the proof of Lemma  \ref{AFC}.
\begin{proof}[Proof of Lemma  \ref{AFC}]
Since every minimal point-distal system is an {\bf AI}-extension of a trivial system, Lemma \ref{AFC} follows from Lemma \ref{admissible ext}, Lemma \ref{EE}, Lemma \ref{AE} and Lemma \ref{inv lim}.
\end{proof}

\section{Unconditonal results}
In this section, we establish Theorem \ref{mainA}, i.e. any totally minimal systems is $\mathbb{P}$-transitive without assuming the Dickson-Hardy-Littlewood conjecture. 

First we outline the ideas of the proof. 
Let $(X,T)$ be a totally minimal system. By Theorem \ref{structure} we have the following diagram:
$$
\xymatrix{
  X          & X_\infty \ar[d]^{\pi_\infty} \ar[l]^{\theta} \\
                & Y_\infty             }
$$
where $X_\infty$ is a proximal extension of $X$ and a RIC
weakly mixing extension of the strictly PI-system $Y_\infty$.

We first deal with PI systems using of quasi-disjointness (Theorem \ref{P-tran for PI}), and prove $Y_\infty$ is $\mathbb{P}$-transitive. Then we deal with RIC weakly mixing extensions. We show that the property $\mathbb{P}$-transitivity can be lifted by RIC weakly mixing extension $\pi_\infty$ (Theorem \ref{P-tran for RIC}). Thus $X_\infty$ is $\mathbb{P}$-transitive. As the factor of $X_\infty$, $X$ is also $\mathbb{P}$-transitive, which completes the proof of Theorem \ref{structure}.
To get better understanding of the proof of Theorem \ref{P-tran for RIC}, we show that any minimal weakly mixing system is $\mathbb{P}$-transitive using Maynard's result (Theorem \ref{prime orbit in wm}).

\subsection{Totally minimal PI systems}

A {\it joining} of two systems $(X,T)$ and $(Y,S)$ is a closed subset of $X\times Y$ that is invariant under $T\times S$ and projects fully onto $X$ and $Y$. The following notion of quasi-disjointness is introduced in \cite{XY 27}.

\begin{defn}[Quasi-disjointness]
Two minimal systems $(X,T)$ and $(Y,S)$ are {\it quasi-disjoint} if $X\times Y$ is the only joining of $X$ and $Y$ that projects onto $X_{eq}\times Y_{eq}$, where $X_{eq}$ and $Y_{eq}$ are the maximal equicontinuous factors of $X$ and $Y$, respectively. 
\end{defn}

If a minimal system is quasi-disjoint itself, then one can deduce the property $A$-transitivity to its maximal equicontinuous factor.

\begin{lem}\label{QD-A-tran}
Let $(X,T)$ be a minimal system satisfying $X$ is quasi-disjoint itself. Then for any subset $A\subset \mathbb{N}$, the following assertions are equivalent:
\begin{enumerate}
\item $(X,T)$ is $A$-transitive;
\item $(X_{eq}, R)$ is $A$-transitive.
\end{enumerate}
\end{lem}
\begin{proof}
It suffices to show that (2) implies (1). Let $\pi: X\rightarrow X_{eq}$ be the factor map and let 
\[ J=\overline{\{(x, T^{n}x):x\in X, n\in A\}}.\]
Then $J$ is a joining of $X$ and $X$. Since $(X_{eq}, R)$ is $A$-transitive, we conclude that $\pi\times\pi(J)=X_{eq}\times X_{eq}$. Since $X$ is quasi-disjoint from itself, we have $J=X\times X$. By Lemma \ref{char of A-tran}, $(X,T)$ is $A$-transitive.
\end{proof}

One of main results in \cite{XY 27} is the following:

\begin{thm}[\cite{XY 27}] \label{QD}
Every minimal PI system is quasi-disjoint from all minimal systems.
\end{thm}

Thus by Lemma \ref{QD-A-tran} and Theorem \ref{QD}, we have the following result.

\begin{thm}\label{P-tran for PI}
Every totally minimal PI system is $\mathbb{P}$-transitive.
\end{thm}
\begin{proof}
Let $(X,T)$ be a totally minimal PI system. By Theorem \ref{QD}, $X$ is quasi-disjoint from itself. Then it follows from Lemma \ref{QD-A-tran} that the $\mathbb{P}$-transitivity of $(X,T)$ is equivalent to the $\mathbb{P}$-transitivity of $(X_{eq},R)$.  Clearly, $(X_{eq},R)$ is $\mathbb{P}$-transitive (see \cite{GT 12} for example).
\end{proof}

\subsection{Minimal weakly mixing systems}\label{MT-WM}

In this subsection, we show that minimal weakly mixing systems are $\mathbb{P}$-transitive. This will help the readers to get better understanding about the ideas in the next subsection.

\begin{thm}\label{prime orbit in wm}
Let $(X,T)$ be a  minimal weakly mixing system. Then there is some point $x\in X$ such that $\overline{orb_{\mathbb{P}}(x)}=X$.
\end{thm}

\begin{proof}
By Lemma \ref{transitive}, it suffices to show that for any nonempty open sets $U, V$ of $X$, $(\bigcup_{p\in\mathbb{P}}T^{-p}U)\cap V\neq\emptyset$.
Since $(X,T)$ is minimal, there is some $k\in \N$ such that
\begin{equation}\label{eq4}
X=T^{-1}V\cup \cdots\cup T^{-k}V.
\end{equation}

Let $m=m_{k}$ be given by Theorem \ref{MT} and let
\[ w= \prod_{p\in\mathbb{P}, p\leq m}p.\]
By Lemma \ref{W-trick}, $\{ wn, 2wn,\ldots, mwn\}$ is admissible,  for any  $n\in \mathbb{N}$.

\medskip
\noindent{\bf Claim}. There are $1\leq a_1<a_2<\cdots<a_k\leq m$ and a piecewise syndetic subset $B\subset \mathbb{N}$ such that for each $n\in B$ there is some $n'\in \N$ such that
\[ \{a_1wn,a_2wn,\ldots,a_kwn\}+n'\subset\mathbb{P}.\]
\begin{proof}[Proof of Claim]
Since $\{ wn, 2wn,\ldots, mwn\}$ is admissible for each $n\in \mathbb{N}$. Then it follows from Theorem \ref{MT} that there are $1\leq a_{1}(n)<a_{2}(n)<\cdots<a_{k}(n)\leq m$ such that $\{ a_{1}(n)wn,a_{2}(n)wn,\ldots,a_{k}(n)wn\}$ can be translated into the primes simultaneously.

This leads to a partition of  $\mathbb{N}$. That is, the partition follows from
\[\phi: \mathbb{N}\rightarrow \{(a_1,\ldots,a_k): 1\leq a_{1}<a_{2}<\cdots<a_{k}\leq m\}, \ \ n\mapsto (a_1(n),\ldots,a_k(n)).\]
That is, $\N=\bigcup_{1\leq a_{1}<a_{2}<\cdots<a_{k}\leq m}\phi^{-1}\big((a_1,\ldots,a_k)\big)$.
Since $\F_{ps}$ has Ramsey property, there are $1\leq a_1<a_2<\cdots<a_k\leq m$ and a piecewise syndetic subset $B\subset \mathbb{N}$ such that for each $n\in B$, $\{ a_{1}wn,a_{2}wn,\ldots,a_{k}wn\}$ can be translated into the primes simultaneously.
\end{proof}

Take $1\leq a_1<a_2<\ldots<a_k\leq m$ be as in the Claim. Now by Lemma \ref{ts for weakly mixing-s}, the set
\[C:= \{n\in\mathbb{Z}: U \cap T^{-a_1wn}T^{-1}U\cap \cdots\cap T^{-a_{k}wn}T^{-k}U\neq\emptyset\}\]
is a thickly syndetic subset of $\mathbb{Z}$. It follows that $C\cap \N$ is a thickly syndetic subset of $\mathbb{N}$, and we have $B\cap C=B\cap (C\cap\N )\neq\emptyset$. Hence there is some $n\in C\cap \N$ such that $(a_1wn,\ldots,a_kwn)$ is prime-producing.

Take $x\in U \cap T^{-a_1wn}T^{-1}U\cap \cdots\cap T^{-a_{k}wn}T^{-k}U$.  Then for any $\ell\in\mathbb{Z}$, there is some $i\in\{1,\ldots, k\}$ such that $T^{\ell}x\in T^{-i}V$. Thus $T^{i}x\in T^{-\ell}V$. But $T^{i}x\in T^{-a_{i}wn}U$. This implies that
\[ (T^{-a_1wn}U\cup T^{-a_{2}wn}U\cup\cdots\cup T^{-a_{k}wn}U)\cap (T^{-\ell}V)\neq\emptyset, \forall \ell\in\mathbb{Z}.\]

Now since $(a_1wn,\ldots,a_kwn)$ is prime-producing, there is some $\ell\in\N$ such that 
\[\{a_1wn,\ldots,a_kwn\}+\ell\subset\mathbb{P}.\]
 Thus we have that
\[ (\bigcup_{p\in\mathbb{P}}T^{-p}U)\cap V\neq\emptyset.\]
This completes the proof.
\end{proof}

\subsection{RIC weakly mixing extensions}
In this subsection, we show that $\mathbb{P}$-transitivity is preserved under taking RIC weakly mixing extensions.

\begin{lem}\label{RIC 1}
Let $\pi: (X,T)\rightarrow (Y,S)$ be an extension between minimal systems. For any nonempty open subsets $U,V$ of $X$, there are finitely many mutually distinct positive integers $a_1,a_2,\ldots,a_r$ and $x\in X$ such that
\[\pi(U)\cap\bigcap_{i=1}^{r}\pi(T^{-a_i}U)\neq\emptyset \ \ \text{ and }\ \  \pi^{-1}\pi(x)\subset \bigcup_{i=1}^{r}T^{-a_i}V .\]
\end{lem}
\begin{proof}
Since $(X,T)$ is minimal, there is some $k\in\mathbb{N}$ such that $V\cap T^{-k}U\neq\emptyset$. Set 
\[ W:=V\cap T^{-k}U.\]
Take $x\in W$ and put $y=\pi(x)$. Then $T^{k}x\in U$. By the minimality of $(X,T)$ and the compactness of $\pi^{-1}(y)$, there are distinct positive integers $a_1,a_2,\ldots,a_r$  such that 
\[ \pi^{-1}(y)\subset \bigcup_{i=1}^{r}T^{-a_i}W\subset\bigcup_{i=1}^{r}T^{-a_i}V, \ \text{and}\ \pi^{-1}(y)\cap T^{-a_i}W \neq\emptyset,\ \forall 1\le i\le r.\]
There are $\{x_1,\ldots,x_r\}\subset \pi^{-1}(y)$ such that $T^{a_i}x_i\in W$ for each $i\in\{1,\ldots,r\}$. Then 
\[ S^{k}y=\pi(T^{k}x_i)\in \pi(T^{k}T^{-a_i}W)\subset \pi(T^{k}T^{-a_i}T^{k}U)=\pi(T^{-a_i}U),\]
for each $i\in\{1,\ldots,r\}$. Therefore,
\[ S^{k}y\in \pi(U)\cap\bigcap_{i=1}^{r}\pi(T^{-a_i}U).\]
Hence $\pi(U)\cap\bigcap_{i=1}^{r}\pi(T^{-a_i}U)\neq\emptyset$.
\end{proof}

Let $(X,T)$ be a system and $k\in\mathbb{N}$. For nonempty open subsets $U_0,U_1,\ldots, U_k$ of $X$ and $\vec{a}=(a_1,\ldots,a_k)\in\mathbb{N}^{k}$, set
\[ N_{\vec{a}}(U_0, U_1,\ldots, U_k):=\{n\in\mathbb{N}:  U_0\cap T^{-a_1n}U_1\cap \cdots\cap T^{-a_kn}U_k\neq\emptyset\}.\]

\begin{lem}\label{RIC 2}
Let $\pi: (X,T)\rightarrow (Y,S)$ be an RIC weakly mixing extension between minimal systems. Then for any $k\in\mathbb{N}$ and $\vec{a}=(a_1,\ldots,a_k)\in\mathbb{N}^{k}$ with $a_1,\ldots, a_k$ pairwise distinct, the collection 
\begin{equation}\label{N_a-filter}
\Big\{N_{\vec{a}}(U_0,U_1,\ldots,U_k): U_0,U_1,\ldots,U_k \text{ are open subsets of } X \text{ and } \bigcap_{i=0}^{k}\pi(U_i)\neq\emptyset\Big\}
\end{equation}
generates a filter. 
\end{lem}
\begin{proof} Recall that $\pi$ is open and $R_{\pi}\subset {\bf RP}^{[\infty]}(X)$. We may assume that $1\leq a_1<a_2<\cdots<a_k$.

\noindent{\bf Claim}. For any open sets $U_0, U_1,\ldots,U_k$ of $X$, if $\cap_{i=0}^{k}\pi(U_i)\neq\emptyset$ then 
\[N_{\vec{a}}(U_0,U_1,\ldots,U_k)\neq\emptyset.\]

\begin{proof}[Proof of Claim]
By Theorem \ref{saturated thm}, there is some point $y\in \cap_{i=0}^{k}\pi(U_i)$
such that $\overline{\mathcal{O}}_{T\times \cdots\times T^{a_k}}(x^{(a_k)})$ is $\pi^{(a_k)}$-saturated. Let 
\[V_i=\begin{cases}
U_0, & i\in\{0,1,2,\ldots, a_k\}\setminus\{a_1,a_2,\ldots,a_k\}\\
U_j, & i=a_j.
\end{cases}\]

Thus there is some $n\in\mathbb{N}$ such that 
\[ V_0\cap T^{-n}V_1\cap T^{-2n} V_2\cap\cdots T^{-a_kn}V_{a_k}\neq\emptyset.\]
In particular, $U_0\cap T^{-a_1n}U_1\cap T^{-a_2n}U_2\cap\cdots\cap T^{-a_kn}U_k\neq\emptyset$ and hence 
\[n\in N_{\vec{a}}(U_0, U_1,\ldots,U_k)\neq\emptyset.\]
The proof of Claim is complete.
\end{proof}

To show the finite intersection property, let $U_0, U_1,\ldots, U_k$ and $V_0,V_1,\ldots, V_k$ be open sets of $X$ satisfying 
\[\bigcap_{i=0}^{k}\pi(U_i)\neq\emptyset\ \ \text{and } \ \  \bigcap_{i=0}^{k}\pi(V_i)\neq\emptyset.\] 
Since $R_{\pi}^{k+1}=\{(x_0,x_1,\ldots,x_k)\in X^{k+1}: \pi(x_0)=\pi(x_1)=\cdots=\pi(x_k)\}$ is transitive under $T\times T\cdots\times T$, there is some $m\in \mathbb{N}$ such that 
\[ R_{\pi}^{k+1}\cap \Big((U_0\cap T^{-m}V_0)\times (U_1\cap T^{-m}V_1)\times\cdots\times (U_k\cap T^{-m}V_k)\Big)\neq\emptyset.\]
Let 
\[ W_i=U_i\cap T^{-m}V_i, i=0,1,\ldots,k.\]
Take $(x_0,x_1,\ldots,x_k)\in R_{\pi}^{k+1}\cap \Big((U_0\cap T^{-m}V_0)\times (U_1\cap T^{-m}V_1)\times\cdots\times (U_k\cap T^{-m}V_k)\Big)$ and put $y=\pi(x_0)=\pi(x_1)=\cdots=\pi(x_k)$. Then $y\in \bigcap_{i=0}^{k}\pi(W_i)\neq\emptyset$ and 
\[N_{\vec{a}}(W_0,W_1,\ldots,W_k)\subset N_{\vec{a}}(U_0,U_1,\ldots,U_k)\cap N_{\vec{a}}(V_0,V_1,\ldots,V_k).\]
This shows the collection (\ref{N_a-filter}) generates a filter.
\end{proof}

\begin{thm}\label{strong MT}
For any $k\in\N$, there is some $m=m_{k}\in\mathbb{N}$ such that for any totally minimal PI system $(X,T)$, any nonempty open subsets $U,V\subset X$, there is a nonempty open subset $U'\subset U$ such that  any admissible sequence $(h_1,h_2,\ldots,h_{m})$, there are $1\leq i_1<\ldots<i_k
\leq m$ such that for any nonempty open subset $U''\subset U'$ there is some   $\ell\in \mathbb{N}$ such  that
\[\{h_{i_1},\ldots,h_{i_k}\}+\ell\subset\mathbb{P} \text{ and } T^{-\ell}U''\cap V\neq\emptyset.\]
\end{thm}

The proof of Theorem \ref{strong MT} is left in the next subsection.

For an admissible set $A\subset \mathbb{Z}$, let 
\[ A_{\mathbb{P}}:=\{ n\in\mathbb{N}: A+n\subset\mathbb{P}\}.\]

\begin{thm}\label{P-tran for RIC}
Let $\pi: (X,T)\rightarrow (Y,S)$ be an RIC weakly mixing extension between totally minimal systems.  If $(Y,S)$ is a totally minimal PI system, then $(X,T)$ is $\mathbb{P}$-transitive.
\end{thm}
\begin{proof}
Fix nonempty open subsets $U,V$ of $X$.  By Lemma \ref{RIC 1}, there are distinct positive integers $b_1,\ldots, b_k$ such that 
\[\pi(U)\cap \bigcap_{i=1}^{k}\pi(T^{-b_i}U)\neq\emptyset \ \ \text{ and }\ \  \pi^{-1}\pi(x)\subset \bigcup_{i=1}^{k}T^{-b_i}V .\]
Recall that $\pi$ is open. The set $\widetilde{U}:=\pi(U)\cap\bigcap_{i=1}^{k}\pi(T^{-b_i}U)$ is a nonempty open subset of $Y$.  Furthermore, take an open neighborhood $\widetilde{V}'$ of $\pi(x)$ such that 
\[ \pi^{-1}(\widetilde{V}') \subset \bigcup_{i=1}^{k}T^{-b_i}V.\]

Let $m=m_k$ be as in Theorem \ref{strong MT}. Applying Theorem \ref{strong MT} with the system $(Y,S)$ and nonempty open sets $\widetilde{U}, \widetilde{V}'$, there is a nonempty open subset  $\widetilde{U}'\subset \widetilde{U}$ such that for any admissible set $\{h_1,\ldots,h_m\}$, there are some $1\leq i_1<\ldots<i_k\leq m$ for any nonempty open subset $\widetilde{U}''\subset \widetilde{U}'$ there is some  $\ell\in\mathbb{N}$ such that
\[\{h_{i_1},\ldots,h_{i_k}\}+\ell\subset \mathbb{P} \ \text{ and } \ S^{-\ell}\widetilde{U}''\cap  \widetilde{V}'\neq\emptyset.\]
 
 Set
\[ w= \prod_{p\in\mathbb{P}, p\leq m}p.\]
It is clear that $\{ wn, 2wn,\ldots, mwn\}$ is admissible,  for any  $n\in \mathbb{N}$. Thus, for any $n\in\mathbb{N}$, there are $1\leq a_{1}(n)<a_{2}(n)<\cdots<a_{k}(n)\leq m$  such that for any nonempty open subset $\widetilde{U}''\subset \widetilde{U}'$ there is some  $\ell_{n,\widetilde{U}''}$ such that 
\[\{ a_{1}(n)wn,a_{2}(n)wn,\ldots,a_{k}(n)wn\}+\ell_{n,\widetilde{U}''}\subset\mathbb{P} \text{ and } S^{-\ell_{n,\widetilde{U}''}}\widetilde{U}'' \cap  \widetilde{V}'\neq\emptyset.\]

This leads to a finite partition of  $\mathbb{N}$. That is, the partition follows from
\[\phi: \mathbb{N}\rightarrow \{(a_1,\ldots,a_k): 1\leq a_{1}<a_{2}<\cdots<a_{k}\leq m\}, \ \ n\mapsto (a_1(n),\ldots,a_k(n)).\]
That is, $\N=\bigcup_{1\leq a_{1}<a_{2}<\cdots<a_{k}\leq m}\phi^{-1}\big((a_1,\ldots,a_k)\big)$.
Set 
\[U^{*}=U\cap \pi^{-1}(\widetilde{U}').\]
Then we have
\[\pi(U^{*})\cap \bigcap_{i=1}^{k}\pi(T^{-b_i}U)\neq\emptyset. \]

For each ${\bf a}=(a_1,\ldots,a_k)$ with $1\le a_1<\cdots<a_k\le m$, let
\[
E_{\bf a}:=N_{(a_1w,\ldots,a_kw)} (U^*,T^{-b_1}U,\ldots,T^{-b_k}U).
\]
Set $W_0=U^*$, $W_{a_i}=T^{-b_i}U$ for $1\leq i\leq k$, and $W_j=X$ at all remaining coordinates. Then
\[
E_{\bf a} =N_{(w,2w,\ldots,mw)}(W_0,\ldots,W_m),\ \  \bigcap_{j=0}^{m}\pi(W_j)=\widetilde{U}'\neq\emptyset.
\]
By Lemma   \ref{RIC 2} applied to $(w,2w,\ldots,mw)$, all the sets $E_{\mathbf a}$ belong to the same filter. Since there are only finitely many choices of ${\bf a}$, we may choose
\[
n_0\in\bigcap_{1\leq a_1<\cdots<a_k\leq m} E_{\bf a}.
\]
Write ${\bf a}=(a_1,\ldots,a_k)=\phi(n_0)$. Then $n_0\in E_{{\bf a}}$, and by the definition of $\phi$, for every
nonempty open subset $\widetilde{U}''\subseteq\widetilde{U}'$ there exists $\ell\in\mathbb N$ such that
\[
\{a_1wn_0,\ldots,a_kwn_0\}+\ell\subseteq\mathbb P \quad\text{and}\quad S^{-\ell}\widetilde U''\cap\widetilde V'\neq\emptyset.
\]

Thus 
\[ W:=U^{*}\cap T^{-a_1wn_0}T^{-b_1}U\cap T^{-a_2wn_0}T^{-b_2}U\cap\cdots\cap T^{-a_kwn_0}T^{-b_k}U\neq\emptyset.\]

Note that $\pi(W)\subset \pi(U^{*})\subset \widetilde{U}'$ is a nonempty open set. Thus there is some $\ell\in\mathbb{N}$ such that
\[\{a_1wn_0, a_2wn_0,\ldots,a_kwn_0\}+\ell\subset\mathbb{P} \text{ and }
S^{-\ell}\pi(W) \cap  \widetilde{V}'\neq\emptyset,\]
 and hence $T^{-\ell}W\cap \pi^{-1}( \widetilde{V}')\neq\emptyset$. In particular,
\[ T^{-\ell}(T^{-a_1wn_0}T^{-b_1}U\cap T^{-a_2wn_0}T^{-b_2}U\cap\cdots\cap T^{-a_kwn_0}T^{-b_k}U) \cap (\bigcup_{i=1}T^{-b_i}V)\neq\emptyset,\]
since $\pi^{-1}(\widetilde{V}') \subset \bigcup_{i=1}^{k}T^{-b_i}V$.
Thus there is some $i\in\{1,\ldots,k\}$ such that 
\[ T^{-\ell}T^{-a_iwn_0}T^{-b_i}U\cap T^{-b_i}V\neq\emptyset.\]
Therefore 
\[ T^{-(a_iwn_0+\ell)}U\cap V\neq\emptyset.\]
Recall that $a_iwn_0+\ell\in\mathbb{P}$. Therefore there is some $p\in\mathbb{P}$ such that $T^{-p}U\cap V\neq\emptyset$. This shows that $(X,T)$ is $\mathbb{P}$-transitive.
\end{proof}

\subsection{Proof of Theorem \ref{strong MT}}
We first give some lemmas. 
\begin{lem}\label{cover+disjoint}
Let $Z$ be a compact metrizable abelian group and $R$ be a minimal rotation on $Z$. For any nonempty open subset $W$ of $X$, there is are nonempty open subsets $U,V$ of $X$ with $V\subset U\subset W$ and finitely many positive integers $n_1,\ldots, n_r$ such that 
\begin{enumerate}
\item $X=\bigcup_{i=1}^{r}R^{-n_i}U$,
\item $R^{-n_1}\overline{V}, \ldots, R^{-n_r}\overline{V}$ are mutually disjoint,
\item $R^{-n_i}\overline{U}\cap R^{-n_j}\overline{V}=\emptyset, \forall i\neq j$.
\end{enumerate}
\end{lem}
\begin{proof}
Let $\rho$ be an invariant metric on $Z$. Take $z_0\in W$ and assume $B(z_0, \epsilon)\subset W$ for some $\epsilon>0$. Let $D=\{x_1,\ldots,x_n\}$ be a maximal $\epsilon$-separated subset of $Z$. That is $\rho(x_i,x_j)\geq \epsilon$ for any $i\neq j$ and $\rho(x, D)< \epsilon$ for any $x\in Z$. This implies that $ Z=\bigcup_{i=1}^{r} B(x_i, \epsilon)$.

Let $\delta=\max_{x\in Z} \rho(x,D)$. Then $\delta<\epsilon$. Take $\eta\in (0,\frac{\epsilon-\delta}{3})$. Since $R$ is a minimal rotation, there are positive integers $n_1,\ldots,n_r$ such that 
\[ \rho(R^{-n_i}z_0, x_i)<\eta, \ \ \forall i\in\{1,\ldots, r\}.\]
For any $x\in X$, there is some $i\in\{1,\ldots,r\}$ such that $\rho(x, x_i)\leq \delta$ and hence 
\[ \rho(x, R^{-n_i}z_0) \leq \delta+\eta.\]
Further, 
\begin{equation}\label{disjoint cover-eq1}
 \rho(R^{-n_i}z_0,R^{-n_j}z_0)\geq \rho(x_i,x_j)-2\eta\geq \epsilon-2\eta,\ \  \forall i\neq j.
 \end{equation}
Noting that $\delta+\eta<\epsilon-2\eta$, we can take $a\in (\delta+\eta, \epsilon-2\eta)$. Set $U=B(z_0, a)$ and we have
\[ Z=\bigcup_{i=1}^{r} B(R^{-n_i}z_0, a)=\bigcup_{i=1}^{r} R^{-n_i}U.\]
Now take $0<b<\min\left\{a, \epsilon-2\eta-a\right\}$ and take $V=B(z_0, b)$.  According to (\ref{disjoint cover-eq1}), we conclude that $R^{-n_1}\overline{V}, \ldots, R^{-n_r}\overline{V}$ are mutually disjoint and  $R^{-n_i}\overline{U}\cap R^{-n_j}\overline{V}=\emptyset, \forall i\neq j$. Clearly, $V\subset U\subset W$. This completes the proof.
\end{proof}

\begin{lem}\label{admissible equi}
Let $(Z,R)$ be a minimal rotation on a compact connected metrizable abelian group $Z$. For any $z_0\in Z$, any nonempty open subsets $V$ of $Z$ such that for any any admissible set $H$ with $|H|=m_k$, there is a subset $A\subset H$ with $|A|=k$,  an infinite subset $C_{A}\subset \mathbb{N}$ and a nonempty open subset $V'\subset V$ such that 
\begin{enumerate}
\item $A+n\subset \mathbb{P}$ for each $n\in C_{A}$ and
\item $ V'\subset \overline{\{ R^{n}z_0: n\in C_{A}\}}$. 
\end{enumerate}
\end{lem}

\begin{proof}
Let $z_0\in Z$ and let $V$ be a nonempty open subset of $Z$. Let $H$ be an admissible set with $|H|$ and set 
\[ F_{H}:=\{n\in\mathbb{N}: |(H+n)\cap \mathbb{P}|\geq k\}.\]
We claim that $\{R^{n}z_0: n\in F_{H}\}$ is dense in $Z$. For this we fix a  nonempty open set $U$ in $Z$.  Since $Z$ is a compact connected metrizable abelian group, there is a homomorphism $\pi: Z\rightarrow \mathbb{T}^{d}$ and a nonempty open subset $W$ of $\mathbb{T}^{d}$ such that $\pi^{-1}(W)\subset U$.  Suppose that $R$ is determined by some $a\in Z$, that is $R(z)=z+a$. Then 
\[R_{\pi(a)}: \mathbb{T}^{d}\rightarrow \mathbb{T}^{d}, x\mapsto x+\pi(a)\]
 is a minimal rotation on $\mathbb{T}^{d}$. By Lemma \ref{torus}, $\{R_{\pi(a)}^{n}\pi(z_0): n\in F_{H}\}$ is dense in $\mathbb{T}^{d}$. Thus there is some $n\in F_{H}$ with $R_{\pi(a)}^{n}\pi(z_0)\in W$. This implies that $R^{n}z_0\in \pi^{-1}(W)\subset U$. The arbitrary choice of $U$ shows our claim.

For each $A\subset H$ with $|A|=k$ and $N\in\mathbb{N}$, set 
\[\Sigma_{A,N}:=\overline{\{R^{n}z_0: n\geq N, n\in A_{\mathbb{P}}\}} \ \ \text{ and } \ \ \Sigma_{A}:=\bigcap_{N=1}^{\infty}\Sigma_{A, N}, \]
where $A_{\mathbb{P}}=\{n\in\mathbb{N}: A+n\subset \mathbb{P}\}$.
Since $F_{H}=\bigcup_{A\subset H, |A|=k} A_{\mathbb{P}}$, for any $N\in\mathbb{N}$ we have 
\begin{equation}\label{eq-equi-0}
  Z=\overline{\{R^{n}z_0: n\in F_{H}\}}=\overline{\{R^{n}z_0: n\in F_{H}, n\geq N\}}=\bigcup_{A\subset H, |A|=k} \Sigma_{A, N}.
\end{equation}
Since $\Sigma_{A,N}$ decreases when $N$ increases, we claim that 
\begin{equation}\label{eq-equi-1}
Z=\bigcup_{A\subset H, |A|=k} \Sigma_{A}.
\end{equation}
Indeed, if there is some $z\in Z$ such that $z\notin\bigcup_{A\subset H, |A|=k} \Sigma_{A}$, then for each $A\subset H, |A|=k$ there is some $N_{A}$ such that $z\notin \Sigma_{A,N_{A}}$. Let $M=\max\{N_{A}: A\subset H, |A|=k\}$. Then $z\notin \bigcup_{A\subset H, |A|=k} \Sigma_{A, M}$, which contradicts (\ref{eq-equi-0}).

Note that  each $\Sigma_{A}$ is closed.  Take a nonempty open subset $W$ with $\overline{W}\subset V$. 
By (\ref{eq-equi-1}), we have $\overline{W}=\bigcup_{A\subset H, |A|=k} \Sigma_{A}\cap \overline{W}$. By Baire Theorem, there is some $A\subset H, |A|=k$ and a nonempty open set $V'\subset V$ such that $V'\subset \Sigma_{A}\cap \overline{W}$. Now let 
\[ C_{A}=\{n\in A_{\mathbb{P}}: R^{n}z_0\in V\}.\]
Clearly,  $A+n\subset\mathbb{P}$ for each $n\in C_{A}$ and 
\[ \Sigma_{A}\cap \overline{W}\subset \bigcap_{N=1}^{\infty}\overline{\{R^{n}z_0: n\in C_A, n\geq N\}}.\]
Thus 
\[V'\subset \bigcap_{N=1}^{\infty}\overline{\{R^{n}z_0: n\in C_A, n\geq N\}}.\]
In particular, $C_{A}$ is infinite and $V'\subset\overline{\{R^{n}z_0: n\in C_A\}}$. The proof is complete.
\end{proof}

\begin{lem}\label{admissible equi-2}
Let $(Z,R)$ be a minimal rotation on a compact connected metrizable abelian group $Z$. For any nonempty open subsets $U, V$ of $Z$ and $z_0\in U$, there is an open neighborhood $U'$ of $z_0$ with $U'\subset U$ such that for any admissible set $H$ with $|H|=m_k$, there is a subset $A\subset H$ with $|A|=k$,  an infinite subset $C_{A}\subset \mathbb{N}$ and a nonempty  open subset $V_{H}'$  such that 
\begin{enumerate}
\item $A+n\subset \mathbb{P}$ for each $n\in C_{A}$ and
\item $ V_{H}'+(z-z_0)\subset V\cap \overline{\{ R^{n}z: n\in C_{A}\}}$ \ for each $z\in U'$. 
\end{enumerate}
\end{lem}

\begin{proof}
Let $U,V$ be any nonempty open subsets of $Z$ and $z_0\in U$.
Take an open neighborhood $U'$ of $z_0$ with $U'\subset U$ and a nonempty open subset $W$ of $Z$ such that 
\[ U'+W\subset V.\]
Now let $H$ be an admissible set with $|H|=m$. By Lemma \ref{admissible equi}, there is some $A\subset H$ with $|A|=k$, an infinite subset $C_{A}\subset \mathbb{N}$ and a   nonempty open subset $W'_{H}$ of $Z$ such that 
\begin{enumerate}
\item $A+C_{A}\subset \mathbb{P}$;
\item $W'_{H}\subset W\cap \overline{\{ R^{n}e_{Z}: n\in C_{A}\}}$,
\end{enumerate}
where  $e_{Z}$ is the unit of the group $Z$.  Now for each $z\in U'$, let
\[ V_{z}=z+W'_{H}.\]
Since $U'+W\subset V$, $V_{z}\subset V$ for any $z\in U'$. Then 
\[ V_{z}\subset V\cap (z+\overline{\{ R^{n}e_{Z}: n\in C_{A}\}})= V\cap \overline{\{ R^{n}z: n\in C_{A}\}}, \ \forall z\in U'. \]
Put $V_{H}'=V_{z_0}$. Then we have $ V_{H}'+(z-z_0)=V_z\subset V\cap \overline{\{ R^{n}z: n\in C_{A}\}}$ for each $z\in U'$. The proof is complete.
\end{proof}

\begin{proof}[Proof of Theorem \ref{strong MT}]
Let $(X,T)$ be a totally minimal PI system and let $\pi: (X,T)\rightarrow (X_{eq}, R)$ be the factor map to its maximal equicontinuous factor.  Then $(X_{eq}, R)$ is also totally minimal and hence $X_{eq}$ is a compact connected metrizable abelian group.

Let $U, V$ be nonempty open subsets of $X$. Since $\pi$ is semiopen, both $\pi(U)$ and $\pi(V)$ has nonempty interiors.  Let $\widetilde{U}$ and $\widetilde{V}$ be the interior of $\pi(U)$ and $\pi(V)$. 

Fix $z_0\in \widetilde{U}$ and let $\widetilde{U}'$ be an open neighborhood of $z_0$ given as in Lemma \ref{admissible equi-2}. That is for any any admissible set $H$ with $|H|=m_k$, there is a subset $A\subset H$ with $|A|=k$,  an infinite subset $C_{A}\subset \mathbb{N}$ and a nonempty  open subset $\widetilde{V}_{H}'$  such that 
\begin{enumerate}
\item $A+n\subset \mathbb{P}$ for each $n\in C_{A}$ and
\item $ \widetilde{V}_{H}'+(z-z_0)\subset \widetilde{V}\cap \overline{\{ R^{-n}z: n\in C_{A}\}}$ for each $z\in \widetilde{U}'$. 
\end{enumerate}

Now let $U'=U\cap \pi^{-1}(\widetilde{U}')$. Let $H$ be an admissible set with $|H|=m_k$ and let $A, C_A, \widetilde{V}_{H}'$ be as above.
 We need to show that for any nonempty open subset $U''\subset U'$ there is some   $\ell\in \mathbb{N}$ such  that
\[A+\ell\subset\mathbb{P} \text{ and } T^{-\ell}U''\cap V\neq\emptyset.\]

Without loss of generality, by Lemma \ref{cover+disjoint}, we may assume that there are positive integers $n_1,\ldots,n_r$  and a nonempty open subset $\widetilde{V}_{H}''\subset \widetilde{V}_H'$ such that 
\begin{enumerate}
\item[(a)] $X_{eq}=\bigcup_{i=1}^{r}R^{-n_i}\widetilde{V}_{H}'$,
\item[(b)] $R^{-n_{1}}\overline{\widetilde{V}_{H}''},\ldots, R^{-n_r}\overline{\widetilde{V}_{H}''}$ are mutually disjoint,
\item[(c)]   $R^{-n_i}\overline{\widetilde{V}_{H}'}\cap R^{-n_j}\overline{\widetilde{V}_{H}''}=\emptyset, \forall i\neq j$.
\end{enumerate}
Otherwise, we can shrink $\widetilde{V}_H'$ a bit. These properties remain hold replacing $\widetilde{V}_{H}',\widetilde{V}_{H}''$ with $\widetilde{V}_{H}'+(z-z_0),\widetilde{V}_{H}''+(z-z_0)$, respectively, for each $z\in X_{eq}$.

Let 
\[ B=\{n\in C_{A}: R^{-n}z_0\in \widetilde{V}_{H}'\}.\]
Then 
\[\widetilde{V}_{H}'\subset\overline{\{ R^{-n}z_0: n\in B\}}.\]
Let
\[ J=\bigcup_{i=1}^{r}(\id\times T^{-n_i})\left(\overline{\{(x, T^{-n}x): x\in X,  n\in B\}}\right)=\overline{\{(x, T^{-n}x): x\in X,  n\in B'\}},\]
where $B'=\bigcup_{i=1}^{r}(n_i+B)$.
Then $(T\times T) (J)=J$ and $J$ is a self-joining of $(X,T)$, i.e. a joining of $(X,T)$ and $(X,T)$.Clearly, 
\[ \bigcup_{i=1}^{r}(\id\times R^{-n_i})\left(\overline{\{(z_0, R^{-n}z_0): n\in B\}}\right)\subset (\pi\times\pi)(J).\]
Since $\widetilde{V}_{H}'\subset\overline{\{ R^{-n}z_0: n\in B\}}$ and $X_{eq}=\bigcup_{i=1}^{r}R^{-n_i}\widetilde{V}_{H}'$, we have 
\[ \bigcup_{i=1}^{r}(\id\times R^{-n_i})\left(\overline{\{(z_0, R^{-n}z_0): n\in B\}}\right)=\{z_0\}\times X_{eq}.\]
Since $(\pi\times\pi)(J)$ is $R\times R$-invariant, we have $(\pi\times\pi)(J)=X_{eq}\times X_{eq}$. By Theorem \ref{QD}, $X$ is quasi-disjoint from itself, and we have $J=X\times X$.   By Lemma \ref{char of A-tran}, $(X,T^{-1})$ is $B'$-transitive. It follows from Lemma that  there is a dense $G_{\delta}$ subset $\Omega$ of $X$ such that 
\begin{equation}\label{dis cover1}
\bigcup_{i=1}^{r}T^{-n_i}\left(\overline{\{ T^{-n}x: n\in B\}}\right)=\overline{\{ T^{-n}x: n\in B'\}}=X,\ \forall x\in \Omega.
\end{equation}

Now let 
\[ F=\{n\in C_{A}: R^{-n}z_0\in \widetilde{V}_{H}''\}.\] 
Then $\overline{\widetilde{V}_{H}''}=\overline{\{ R^{-n}z_0: n\in F\}}$. For any $ z\in \pi(U'') \subset \widetilde{U}'$, one has
\begin{align*}
\overline{\widetilde{V}_{H}''+(z-z_0)}&=\overline{\widetilde{V}_{H}''}+(z-z_0)=\overline{\{ R^{-n}z_0: n\in F\}}+(z-z_0)\\
&=\overline{\{ R^{-n}z_0: n\in F\}+(z-z_0)}\\
&=\overline{\{ R^{-n}z: n\in F\}}.
\end{align*}
 Take $x\in U''\cap \Omega$ and let $z=\pi(x)$. Then we have 
\[\overline{\{ T^{-n}x: n\in F\}}\subset \pi^{-1}\left(\overline{\widetilde{V}_{H}''+(z-z_0)}\right).\]
By the definitions of $B$ and $F$, we have 
 \[\{R^{-n}z_0: n\in B\setminus F\} \subset  \widetilde{V}_{H}'\setminus \widetilde{V}_{H}''\]
 and hence 
 \[\{R^{-n}z: n\in B\setminus F\} \subset \left( \widetilde{V}_{H}'+(z-z_0)\right)\setminus \left(\widetilde{V}_{H}''+(z-z_0)\right).\]
 Thus
 \[\{T^{-n}x: n\in B\setminus F\} \subset \pi^{-1}\left( \widetilde{V}_{H}'+(z-z_0)\right)\setminus \pi^{-1}\left(\widetilde{V}_{H}''+(z-z_0)\right)\]
 and 
 \[\bigcup_{i=1}^{r}T^{-n_i}\left(\overline{\{T^{-n}x: n\in B\setminus F\}}\right) \subset\bigcup_{i=1}^{r}T^{-n_i}\left(\overline{\pi^{-1}\left( \widetilde{V}_{H}'+(z-z_0)\right)\setminus \pi^{-1}\left(\widetilde{V}_{H}''+(z-z_0)\right)}\right).\]

 Write
\[D=\widetilde V_H'+(z-z_0), \qquad E=\widetilde V_H''+(z-z_0).\]
Since $E$ is open, the preceding inclusion gives
\[\overline{\{T^{-n}x:n\in B\setminus F\}} \subseteq \pi^{-1}(\overline D\setminus E). \]
By (c), $R^{-n_i}\overline D\cap R^{-n_j}E=\emptyset$ for $i\neq j$. For $i=j$, the set $R^{-n_i}(\overline D\setminus E)$ is disjoint from $R^{-n_i}E$. It follows that
\[ \left(\bigcup_{i=1}^{r}T^{-n_i} \overline{\{T^{-n}x:n\in B\setminus F\}}\right)
\cap \left(\bigcup_{j=1}^{r}T^{-n_j}\pi^{-1}(E)\right) =\emptyset. \]
Combining this with (\ref{dis cover1}), we obtain
\[ \bigcup_{i=1}^{r}T^{-n_i}\overline{\{T^{-n}x:n\in F\}} \supseteq \bigcup_{j=1}^{r}T^{-n_j}\pi^{-1}(E). \]
Moreover, $\overline{\{T^{-n}x:n\in F\}}\subseteq \pi^{-1}(\overline E)$. Since $E\subseteq D$, condition (c) implies that $T^{-n_1}\pi^{-1}(E)$ is disjoint from $T^{-n_i}\pi^{-1}(\overline E)$ for every $i\ne 1$.
Therefore
\[ \pi^{-1}(E)\subseteq \overline{\{T^{-n}x:n\in F\}}. \]
 
 
 
Since $V\cap \pi^{-1}(E)$ is a nonempty open set, there is some $\ell\in F$ such that $T^{-\ell}x\in V\cap \pi^{-1}(E)$.  Thus we have $A+\ell\subset \mathbb{P}$ and $T^{-\ell}U''\cap V\neq\emptyset.$
This completes the proof of Theorem \ref{strong MT}.
\end{proof}

\subsection{Proof of Theorem \ref{mainA}}
Now we can give the proof of Theorem \ref{mainA}.

\begin{proof}[Proof of Theorem \ref{mainA}]
Let $(X,T)$ be a totally minimal system. Then there is a proximal extension $\pi: (\tilde{X}, \tilde{T})\rightarrow (X,T)$ such that $(\tilde{X}, \tilde{T})$ is  a RIC weakly mixing of a strictly PI system $(Y,S)$. By Lemma \ref{total min for prox ext}, $(\tilde{X}, \tilde{T})$ is a totally minimal system. 
By Theorem \ref{P-tran for PI}, $(Y,S)$ is $\mathbb{P}$-transitive. Further, it follows from Theorem \ref{P-tran for RIC} that $(\tilde{X}, \tilde{T})$ is $\mathbb{P}$-transitive. Passing to the factor, we conclude that $(X,T)$ is $\mathbb{P}$-transitive.
\end{proof}

\appendix 
\section{Recurrence in minimal weakly mixing systems}\label{rec in wm}
In this section, we show a result on the multiple hitting time in minimal weakly mixing systems along the line given in \cite{HSY 19}.  In this section, when $G$ is a group, $(X,G)$ is a $G$-action system; and when $T$ is a homeomorphism, $(X,T)$ is the corresponding $\Z$-action system. 

\begin{lem}\label{simu for abelian}
Let $G$ be an abelian group and $H$ be a subgroup of $G$. Suppose that  $(X,H)$ is weakly mixing. Let $V_1,\ldots, V_k$ be nonempty open set in $X$. Then for any finite subset $\{g_1,\ldots, g_r\}\subset G$, there are nonempty open set $W_1,\ldots, W_k$  and $\{h_1,\ldots,h_r\}\subset H$ such that
\[ h_ig_i W_j\subset V_j, \ \ \forall (i,j)\in\{1,\ldots,r\}\times\{1,\ldots,k\}.\]
\end{lem}

\begin{proof}
We show the lemma by induction on $r$.

For $r=1$, we take $h_1=e$ and $W_j=g_{1}^{-1}V_j$ for each $j=1,\ldots, k$. Then it is clear that $h_1g_1W_j=V_j$ for each $j=1,\ldots,k$.

Now suppose that there are nonempty open sets $U_1,\ldots, U_k$ and $h_1,\ldots, h_{r-1}\in H$ such that
\[ h_ig_i U_j\subset V_j, \ \ \forall (i,j)\in\{1,\ldots,r-1\}\times\{1,\ldots,k\}.\]
Since $(X,H)$ is weakly mixing, $(X^{k}, H)$ is transitive under the diagonal action. Thus there is some $h_r\in H$ such that
\[ (h_rg_{r}U_1\times h_rg_r U_2\times \cdots\times h_rg_r U_k)\cap (V_1\times V_2\times\cdots\times V_k)\neq\emptyset.\]
Now let $W_j=U_i\cap(h_r^{-1}g_{r}^{-1}V_j)$ for each $j=1,\ldots,k$. Then each $W_j$ is a nonempty open set. Since $W_j\subset U_j$ for each $j=1,\ldots,k$, one has
\[ h_ig_i W_j\subset h_ig_i U_i\subset V_j, \ \ \forall (i,j)\in\{1,\ldots,r-1\}\times\{1,\ldots,k\}.\]
While it is clear that $h_rg_r W_j\subset V_j$ for each $j=1,\ldots,k$.  Thus we have shown that
\[ h_ig_i W_j\subset V_j, \ \ \forall (i,j)\in\{1,\ldots,r\}\times\{1,\ldots,k\}.\]
The proof is complete.
\end{proof}

\begin{lem}\label{syn in wm2}
Let $\Gamma$ be an abelian group, and $S\in \Gamma$. Suppose that $(Y,\Gamma)$ is minimal and $(Y,S)$ is weakly mixing. Then for any nonempty open sets $V_0,V_1,\ldots,V_{k}$ and pairwise distinct non-zero integers $a_1,\ldots,a_k$, the set
\[E:=\{n\in\Z: V_0\cap S^{-a_1n}V_1\cap \cdots\cap S^{-a_kn}V_k\neq\emptyset\}\]
is  a syndetic subset of $\Z$.
\end{lem}
\begin{proof}
We show the lemma by induction on $k$.
Since $(Y,\Gamma)$ is minimal and $\Gamma$ is abelian, we conclude that $(Y, S)$ has dense minimal points. Then one has $\{n\in\Z: V_0\cap S^{-a_1n}V_1\neq\emptyset\}$ is syndetic and the case of $k=1$ is followed.

Since $(Y, \Gamma)$ is minimal, there is a finite subset $\{g_1,\ldots,g_r\}\subset \Gamma$ such that
\[Y=g_1^{-1}V_0\cup g_2^{-1} V_0\cup \cdots\cup g_{r}^{-1}V_0.\]

Let $\tau:=S^{a_1}\times S^{a_2}\times\cdots\times S^{a_k}$. For each $g\in \Gamma$, let $g^{(k)}=g\times g\times\cdots\times g$ be the homeomorphism on $Y^{k}$. Let $H=\langle\tau\rangle$ and $G=\langle \tau, \{g^{(k)}: g\in \Gamma\}\rangle$. Then $G$ is abelian. Applying Lemma \ref{simu for abelian}  to $Y^{k}, G, H$, there is a nonempty open set $U_1\times \cdots\times U_k$  and $n_1,\ldots,n_r\in\mathbb{Z}$ such that
\[ \tau^{n_i} g_{i}^{(k)} (U_1\times \cdots\times U_k)\subset V_1\times V_2\times\cdots\times V_k,\ \ \forall i\in\{1,\ldots,r\}.\]

Since $a_2-a_1,\ldots,a_k-a_1$ are distinct nonzero integers, it follows from the induction hypothesis that the set
\[ F:=\{ n\in\Z: U_1\cap S^{-(a_2-a_1)n}U_2\cap\cdots\cap S^{-(a_k-a_1)n} U_{k}\neq\emptyset\}\]
is syndetic.

Now take $m\in F$. Then there is some $x_m\in U_1$ such that $S^{(a_j-a_1)m}x_m\in U_j$ for each $j=1,2,\ldots,k$. Let $y_m=S^{-a_1m} x_m$. There is some $l_{m}\in\{1,\ldots,r\}$ such that $g_{l_m}^{-1}(z_m)=y_m$ with $z_m\in V_0$. Now for each $i=1,\ldots,k$,
\begin{align*}
S^{a_i(m+n_{l_m})} z_m=S^{a_i(m+n_{l_m})} g_{l_m}y_m=S^{a_in_{l_m}}g_{l}S^{(a_i-a_1)m}x_m\in S^{a_in_{l_m}}g_{l}U_i\subset V_i.
\end{align*}
By letting $n=m+n_{l_m}$, we have
\[ z_m\in V_0\cap S^{-a_1n}V_1\cap \cdots\cap S^{-a_kn}V_k.\]
Thus $E\supset \{m+n_{l_m}:m\in F\}$ is syndetic.
\end{proof}

\begin{lem}\label{ts for weakly mixing-s}
Let $(X,T)$ be a minimal weakly mixing system and $a_1,\ldots,a_m$ be distinct nonzero integers. Then for any nonempty open sets $U_0,U_1,\ldots, U_{m}\subset X$, the set
\[ \{n\in \Z: U_0\cap T^{-a_1n}U_1\cap \cdots\cap T^{-a_{m}n}U_{m}\neq\emptyset\}\]
is a thickly syndetic subset of $\Z$.
\end{lem}
\begin{proof}
For each $r\in\Z$, set
\[A_{r}=\{n\in Z: U_0\cap T^{-a_1n}T^{-a_1r}U_1\cap \cdots\cap T^{-a_{m}n}T^{-a_mr}U_{m}\neq\emptyset\}.\]
 We need to show that $\bigcap_{r=0}^{k-1} A_{r}$ is syndetic for any $k\in\N$. As a consequence, this implies our assertion.

 Now we fix $k\in \N$ and let
\begin{align*}
V_0&=U_0\times U_0\times\cdots\times U_0 \ \ (k-\text{times}),\\
V_1&= U_1\times T^{-a_1}U_1\times \cdots\times T^{-(k-1)a_1}U_1,\\
&\cdots\\
V_m&= U_m\times T^{-a_m}U_m\times \cdots\times T^{-(k-1)a_m}U_m.
\end{align*}

Let $S=T\times T\times\cdots\times T (k\text{ times })$. Note that the system $(X^{k}, S)$ is weakly mixing and the system $(X^{k}=N_d(X),\langle S, T\times T^2\times\cdots\times T^{k}\rangle)$ is minimal by Lemma \ref{Glasner}. By Lemma \ref{syn in wm2}, the set
\[ E=\{n\in\Z: V_0\cap S^{-a_1n}V_1\cap\cdots\cap S^{-a_mn}V_m\neq\emptyset\}\]
is syndetic in $\Z$. Clearly, $E\subset \bigcap_{r=0}^{k-1}A_{r}$ and hence we conclude the proof.
\end{proof}


\begin{thebibliography}{SSS}

\bibitem{AMWW 84} J. Auslander, D. McMahon, J. van der Woude,  and T. S. Wu, {\it Weak disjointness and the equicontinuous structure relation}. Ergodic Theory Dynam. Systems 4 (1984), no. 3, 323--351.

\bibitem{BST23} P. Y. Bienvenu, X. Shao and J. Ter\"{a}v\"{a}inen, {\it A transference principle for systems of linear equations, and applications to almost twin prime}. Algebra and Number Theory, {\bf 17: 2} (2023).  

\bibitem{Bour 89} J. Bourgain, {\it Pointwise ergodic theorems for arithmetic sets} (with an appendix by H. Furstenberg, Y. Katznelson and D. Ornstein). Publ. Math. IH\'{E}S.  69 (1989), 5--45.


\bibitem{CKN 07} G. Cairns, A. Kolganova and A. Nielsen, {\it Topological transitivity and mixing notions for group actions}. Rocky Mountain J. Math. 37 (2007), no. 2, 371--397.

\bibitem{Dick 04} L. E. Dickson, {\it A new extension of Dirichlet’s theorem on prime numbers}, Messenger of Math. 33 (1904), 155--161.


\bibitem{EGS 75}R. Ellis, S. Glasner, Shmuel and L. Shapiro, {\it Proximal-isometric flows}. Advances in Math. 17 (1975), no. 3, 213--260.

\bibitem{Fur 63}H. Furstenberg. {\it The structure of distal flows}.  Amer. J. Math., 85 (1963), 477--515. 


\bibitem{Fur81} H. Furstenberg, {\it Recurrence in Ergodic Theory and Combinatorial Number Theory}. M. B. Porter Lectures. Princeton, NJ: Princeton University Press, 1981.

\bibitem{G75} E. Glasner, {\it Proximal flows}.  Lecture Notes in Mathematics, Vol. 517. Berlin: Springer-Verlag, 1976.

\bibitem{Glasner 94} E. Glasner, {\it Topological ergodic decompositions and applications to products of powers of a minimal transformation}. J. Anal. Math. 64 (1994), 241--262.

\bibitem{G03}E. Glasner, {\it Ergodic theory via joinings}. Mathematical Surveys and Monographs, Vol. 101. Providence, RI: American Mathematical Society, 2003

\bibitem{G05}E. Glasner, {\it Topological weak mixing and quasi-Bohr systems}.  Isr. J. Math., 148 (2005), 277--304.

\bibitem{GHSWY 25} E. Glasner, W. Huang, S. Shao, B. Weiss and X. Ye, {\it Topological characteristic factors and nilsystems}. J. Eur. Math. Soc. 27 (2025), no. 1, 279--331.

\bibitem{GT 12} B. J. Green and T. Tao, {\it The M\"{o}bius function is strongly orthogonal to nilsequences}. Ann. of Math. 175 (2012), 541--566.

\bibitem{HL 23} G. H. Hardy and J. E. Littlewood, {\it Some problems of ‘Partitio Numerorum.’ III. On the expression of a number as a sum of primes}, Acta Math. 44 (1923), 1--70.

\bibitem{HKM}B. Host, B.  Kra and A. Maass, {\it Nilsequences and a structure theorem for topological dynamical systems}. Adv. Math. 224 (2010), 103--129.

\bibitem{HSY 19} W. Huang, S. Shao and X. Ye, {\it Topological correspondence of multiple ergodic averages of nilpotent group actions}. J. Anal. Math. 138(2019), No. 2, 687--715. 


\bibitem{KLR 24} A. Kanigowski, M. Lema\'{n}czyk and M. Radziwi\l \l, {\it Prime numver theorem for analytic skew products}. Ann. of Math. (2) 199 (2024), no. 2, 591--705.


\bibitem{Maynard 15} J. Maynard, {\it Small gaps between primes}. Ann. of Math. (2) 181 (2015), no. 1, 383--413.

\bibitem{Mc76} 
D. C. McMahon, \emph{Weak mixing and a note on a
structure theorem for minimal transformation groups}, Illinois J.
Math. {\bf 20} (1976), no. 2, 186--197.




\bibitem{Pavlov 08} R. Pavlov, {\it Some counterexamples in topological dynamics}. Ergod. Theory Dynam. Sys. 28 (2008), 1291--1322.
 
\bibitem{Qiu 23} J. Qiu, {\it Polynomial orbits in totally minimal systems}.  Adv. Math. {\bf 432} (2023), Paper No. 109260.


\bibitem{SU 15} P. Sarnak and A. Ubis, {\it The horocycle flow at prime times}.  J. Math. Pures Appl. (9) 103 (2015), no. 2, 575--618.

\bibitem{SY 12} S. Shao and X. Ye, {\it Regionally proximal relation of order d is an equivalence one for minimal systems and a combinatorial consequence}. Adv. Math. 231 (2012), 1786--1817.

\bibitem{Veech 70}W. A. Veech,  {\it Point-distal flows}. Amer. J. Math. 92 (1970), 205--242.


\bibitem{V77} W. A. Veech, \emph{Topological dynamics}, Bull. Amer. Math. Soc.,
{\bf 83}(1977), 775--830.


\bibitem{V37} I. M. Vinogradov,  {\it The representation of an odd number as a sum of three prime numbers}.  Dokl.
Akad. Nauk SSSR (2) 15(1937), 291--294.

\bibitem{Vries 93} J. de Vries, {\it Elements of topological dynamics}. Mathematics and Its Applications, Vol. 257. Dordrecht: Kluwer Academic Publishers, 1993.

\bibitem{XY 27} H. Xu and X. Ye, {\it Quasi-disjointness in topological dynamics}. J. Differential Equations 486 (2027), Paper No. 114760.

\bibitem{Zhang 14} Y. Zhang, {\it Bounded gaps between primes}, Ann. of Math. 179 (2014), 1121--1174.


\end{thebibliography}
\end{document}